\documentclass{article}[12pt]

\usepackage{xypic,graphicx}
\usepackage{amsmath,amssymb}
\usepackage{theorem}
\usepackage{url}
\newcommand{\sil}[1]{}
\newcommand{\Ref}[1]{(\ref{#1})}
\newcommand{\ds}{\displaystyle}

\newcommand{\halala}[1]{\ensuremath{\mathbb{#1}}}
\newcommand{\Ni}{\halala{N}}
\newcommand{\Ri}{\halala{R}}
\newcommand{\Qi}{\halala{Q}}

\newcommand{\K}{\ensuremath{\mathbb{K} } }
\newcommand{\F}{\ensuremath{\mathbb{F} } }

\newcommand{\q}{\ensuremath{\! : \!}}

\newcommand{\gid}{\ensuremath{G}} 
\newcommand{\gva}{\ensuremath{\mathcal{G}}} 
\newcommand{\sid}[1][]{\ensuremath{P_{#1}}} 
\newcommand{\sva}{\ensuremath{\mathcal{P}}} 
\newcommand{\rhid}{\ensuremath{J}} 
\newcommand{\mfid}{\ensuremath{I}} 
\newcommand{\mfmid}{\ensuremath{M}} 
\newcommand{\oid}{\ensuremath{O}} 
\newcommand{\ova}{\ensuremath{\mathcal{O}}} 

\newcommand{\zva}[1][\infty]{\ensuremath{\mathcal{Z}}}
\newcommand{\dimz}{n}
\newcommand{\dimo}{s}
\newcommand{\dv}{\zeta}
\newcommand{\locs}{\ensuremath{\K[z]_{\sva}}}
\newcommand{\ratif}{\ensuremath{{\K(z)^G}} } 
\newcommand{\algif}{\ensuremath{\overline{\K(z)}^G} } 

\newcommand{\ninv}{\ensuremath{\xi} }
\newcommand{\bz}{\ensuremath{\bar{z}} }

\newcommand{\bl}{\ensuremath{\bar{\lambda}} }
\newcommand{\bm}{\ensuremath{\bar{\mu}} }
\newcommand{\brf}{\ensuremath{{f}} }
\newcommand{\prel}{\ensuremath{\theta} }

\newcommand{\ncs}{\ensuremath{\mathcal{U}}}
\newcommand{\eez}{\exp(\varepsilon v, {\bz})}
\newcommand{\bv}[2]{\frac\partial{\partial #1_{#2}}}
\newcommand{\smr}{f}
\newcommand{\smz}{F}
\newcommand{\Ratif}{\ensuremath{{\Ri(z)^G}} }

\newcommand{\inv}{\ensuremath{\iota}}
\newcommand{\isv}{\ensuremath{ \phi_{\ninv}}}
\newcommand{\icv}{\ensuremath{\bar{\iota}}}
\newcommand{\pr}{\pi}

\newenvironment{proof}{\textsc{proof:}}{\ensuremath{\Box}}

\theoremheaderfont{\scshape}
\theorembodyfont{\slshape}
\newtheorem{exy}{Example}[section]
\newtheorem{propy}[exy]{Proposition}
\newtheorem{coly}[exy]{Corollary}
\newtheorem{theoy}[exy]{Theorem}
\newtheorem{lemy}[exy]{Lemma}
\newtheorem{defy}[exy]{Definition}
\newtheorem{hypy}[exy]{Asumption}

\renewcommand{\sec}[1]{Section~\ref{#1}}
\newcommand{\rprop}[1]{Proposition~\ref{#1}}

\newcommand{\theo}[1]{Theorem~\ref{#1}}
\newcommand{\rlem}[1]{Lemma~\ref{#1}}
\newcommand{\ex}[1]{Example~\ref{#1}}

\newcommand{\hyp}[1]{Asumption~\ref{#1}}

\title{Rational, Replacement and Local  Invariants
of a Group Action}

\author{ \parbox[t]{0.5\textwidth}{\begin{center} Evelyne Hubert \\
INRIA Sophia Antipolis \\ {\small
\url{www.inria.fr/cafe/Evelyne.Hubert}}
          \end{center} }
   \parbox[t]{0.5\textwidth}{\begin{center} Irina A. Kogan \\ North
   Carolina State University \\ {\small
   \url{www.math.ncsu.edu/~iakogan}} \end{center}} }

\begin{document}

\maketitle

\begin{abstract}
The paper presents a new algorithmic construction of  a finite
generating set of rational invariants for  the rational action of 
an algebraic group on the affine space. 
The construction  provides an algebraic counterpart  
of the  moving frame method in differential geometry.
The generating set of rational invariants appear as 
the coefficients of a Gr\"obner basis,  reduction with respect to 
which   allows to express a rational invariant in terms of the
generators. 
The replacement invariants, introduced in the paper, are tuples of algebraic
functions of rational invariants. 
Any invariant, whether rational, algebraic or local,  
can be rewritten in terms
of a replacement invariant by a simple substitution.
\end{abstract}

{\bf Key words:}
 rational and algebraic invariants, algebraic and Lie group actions,
cross-section,  Gr\"obner basis,  moving frame method,
smooth and differential invariatns.
 

\section{Introduction}

We present algebraic constructions for invariants of a rational 
group action on an affine space, and relate them to their counterparts in 
differential geometry. 
The constructions are algorithmic and can  easily be implemented in 
general purpose computer algebra systems or software 
specialized in Gr\"obner basis computations. This is illustrated 
by the \textsc{maple} worksheet available at 
\url{http://www.inria.fr/cafe/Evelyne.Hubert/Publi/rrl_invariants.html}
where the examples of the paper are treated.

The first construction is for the computation of  a generating set 
of rational invariants. This generating set is endowed with a 
simple algorithm to express any rational invariant in terms of them.
The construction comes into two variants.
In the first one we consider the ideal of the graph of the action
as did Rosenlicht 
\cite{rosenlicht56}, Vinberg \& Popov \cite{vinberg89}\footnote{
 We are  indebted to  a referee of the MEGA conference for 
pointing out this reference that motivated us to push further some
of the results.}
, and 
 M{\"u}ller-Quade \& Beth \cite{beth99}\footnote{
 We would like to thank H.~Derksen for suggesting comparison with this
 reference after we made public a first  preprint. }.
We point out the connections with these previous works in the text.
Our proofs are independent and provide an original approach.
We show that the coefficients of a reduced Gr\"obner basis of the
ideal of the graph of the action are invariant.
We prove that these coefficients generate the field of rational 
invariants by exhibiting an
algorithm for rewriting any rational invariant in  terms of them.
The second variant provides a purely algebraic 
formulation  of the  geometric  construction of a
\emph{fundamental set of local  invariants}  on a smooth
manifold proposed by Fels and Olver \cite{olver99}, as a generalization of
Cartan's moving frame method.   
It is also computationally more effective  as we reduce to zero 
the dimension of the polynomial ideal for  which  a reduced
Gr\"obner basis is computed.
This is achieved by adding the ideal of a cross-section  to the ideal
of the graph.

That latter construction allows to 
introduce  \emph{replacement invariants},
the algebraic counterpart of   \emph{normalized invariants}
appearing in the geometric construction. 
A replacement invariant is a tuple of algebraic of
functions of rational invariants.
Any invariant can be trivially rewritten in their terms  by 
substituting the coordinate functions by 
the corresponding invariants from this tuple. 
An  \emph{invariantization} map, a computable isomorphism from the set
of algebraic functions on the  
cross-section to the set of algebraic invariants, is defined in terms
of replacement invariants. 

We use invariantization process  to make
explicit the connection between the present algebraic construction 
and the geometric construction of Fels and Olver  \cite{olver99}. 
We introduce an alternative definition  of smooth invariantization
which, on one hand, generalizes the one given in \cite{olver99} 
and, on the other hand, matches the algebraic construction. 
We  thus provide a bridge between the theory of rational
and algebraic invariants \cite{vinberg89,derksen02} and the
theory of smooth local invariants in differential geometry. 

Diverse fields of application of algebraic invariant theory are
presented in 
\cite[Chapter 5]{derksen02}. Some of
the applications  can be addressed with rational invariants. Their
present construction  together with
the simple rewriting algorithm can bring computational benefits.
An application of the  moving frame method to classical invariant
theory \cite{HilbertEng,Gur64,Sturmfels93} was  proposed in
\cite{O99,Kthesis, BO00, KM02}.  In
these works,  
however,  the geometric formulation of the method is used without adapting
it to the algebraic nature of the problem.  A purely algebraic
formulation of the moving frame method opens new possibilities of its
application in classical invariant theory.

The present algebraic formulation provides a new tool  for the 
investigation of the differential invariants of Lie group actions 
and their applications to differential systems in the line of
\cite{olver99,hubert05,ko03,mansfield01}.  
This larger project motivates our choice
to consider rational actions. Even if we start with an affine or even
linear action on the zeroth order jet space, the prolongation of the
action to the higher order jet spaces is usually rational.

The paper is structured as follows.  In \sec{groupaction} we introduce
the action of an algebraic group on the affine space and the graph of
the action. This leads to a first  construction of a set of generating
rational invariants. 
A second version of the construction is given after the introduction
of the cross-section to the orbits in  \sec{crossection}.
This second construction gives rise to the replacement invariants  in 
\sec{invariantization}, which  are used to  define a 
computable invariantization map.
In \sec{felsolver} we present a geometric construction of local smooth invariants that  generalizes the construction of \cite{olver99} and 
explicitly relates it to the algebraic construction of the previous sections. 
 \sec{examples} provides additional examples.

\textsc{Acknowledgments:}
{We would like to thank  Liz Mansfield, Peter Olver and Agnes
  Szanto  for discussing the ideas of the paper during the workshop
  "Differential Algebra and Symbolic
  Computation"  in Raleigh, April 2004, 
  sponsored in part by NSF grants CCR-0306406 and CCF-0347506. 
  We are grateful to
  Michael Singer for continuing discussion of the project and a number
  of valuable suggestions.} 
\section{Graph of a group action and rational invariants}
\label{groupaction} 

We give a definition of a rational action  of an algebraic group over
a field $\K$ on an affine space, and formulate  two additional
hypotheses necessary to our construction.  
We recall the definition for the graph of the action. It plays a
central  role in  our constructions.  The first variant of the
algorithm for constructing  a generating set of rational invariants, 
together with an algorithm for expressing any rational invariant in 
 terms of them, is presented in this section. 

For exposition convenience we assume that the field $\K$  is
algebraically closed.   
The construction proposed in this section 
relies only on  Gr\"obner basis computations and thus can 
be performed in the field of definition of the data (usually $\Qi$ or
$\F_p$). 
Outside of  \sec{felsolver} the terms
\emph{open}, \emph{close} and \emph{closure} refer to the Zariski
topology.

\subsection{Rational action of an algebraic group} \label{agroupaction:def}

We consider an algebraic group that is defined  as  an algebraic variety $\gva$ in the affine space $\K^l$.
The group operation and the inverse are given by polynomial maps.
The neutral element is denoted by $e$.
We shall consider an action of $\gva$ on an affine
space $\zva=\K^\dimz$. 

Throughout the paper $\lambda=(\lambda_1, \ldots,\lambda_{l})$ and 
$z=(z_1,\ldots, z_{\dimz})$ denote indeterminates while
$\bl=(\bar{\lambda}_1, \ldots, \bar{\lambda}_{l})$ and
$\bz=(\bar{z}_1,\ldots, \bar{z}_{\dimz})$ denote points in
$\gva\subset \K^{l}$ and $\zva=\K^{\dimz}$ respectively.
The coordinate ring  of $\zva$ and $\gva$ are respectively
$\K[z_1, \ldots, z_\dimz]$ and  $\K[\lambda_1, \ldots,\lambda_{l}]/G$
where $G$ is  a radical unmixed dimensional ideal.
By $\bl\cdot\bm$ we denote the image of $(\bl,\bm)$ under the group
operation while $\bl^{-1}$ 
denotes the image of $\bl$ under the inversion map.

\begin{defy} \label{action:def} 
A rational action of an algebraic group $\gva$ on the
  affine space $\zva$ is a rational map $g\colon\gva \times \zva
  \rightarrow \zva $ 
that satisfies the following two properties
\begin{enumerate}
\item $g(e,\bz)=\bz$, $ \forall \bz\in \zva$
 \item \(g(\bm, g(\bl, z))=g(\bm\cdot\bl, z)\), whenever both
 $(\bl,\bz)$ and $(\bm\cdot\bl,\bz)$ are in the domain of definition of
 $g$.
\end{enumerate}
\end{defy}

A rational action is thus uniquely determined by a $\dimz$-tuple of
rational functions of $\K(\lambda,z)$ whose domain of definition 
is a dense open set of $\gva\times\zva$. We can bring these rational
functions to their least common denominator $h\in \K[\lambda,z]$ 
without affecting the domain of 
definition. In the rest of the paper the action is thus given by
\begin{equation}\label{action} g(\bl,\bz)=\left(
{g_1(\bar{\lambda}, \bar{z})}, \ldots, {g_{\dimz}(\bar{\lambda},
\bar{z})} \right)\mbox{ for } g_1, \ldots, g_{\dimz} \in
 h^{-1}\K[\lambda_1, \ldots,\lambda_{l},z_1,\ldots,z_{\dimz}]
\end{equation}

\begin{hypy} \label{groupaction:hyp} 
We make the additional assumptions \begin{enumerate}
\item for all $\bz\in \zva$, $h(\lambda,\bz)\in \K[\lambda]$ 
  is not a zero-divisor of $G$. 
  This says that the domain of definition of 
  $g_{\bz} \colon \bl \mapsto g(\bl,\bz)$ contains a non-empty open set of each
  component of $\gva$. 
\item for all $\bl\in \zva$, $h(\bl,z)\in \K[z]$ is different 
  from zero. In other words, for every element  
  $\bl\in\gva$ there exists $ \bz\in\zva$, such that 
  $(\bl,\bz)$ is in the domain of definition $g$. 
\end{enumerate}
\end{hypy}

The following three examples  serve as illustration throughout the text.

\begin{exy} \label{scaling:def} \textsc{Scaling}.
 Consider the multiplicative group given by
$\gid= (1-\lambda_1\lambda_2)\subset \K[\lambda_1, \lambda_2]$.
The neutral element is $(1,1)$ and
$(\bm_1,\bm_2)\cdot (\bl_1,\bl_2)^{-1}=(\bm_1\bl_2,\bm_2\bl_1)$.
We consider the scaling action of this group on $\K^2$.
It is given by the following polynomials of
$\K[\lambda_1,\lambda_2,z_1,z_2]$:
$g_1 = \lambda_1 z_1, \quad g_2=\lambda_1 z_2.$
\end{exy}

\begin{exy} \label{translation:def} \textsc{translation+reflection.}
Consider the group that is the cross product of the additive group and
the group of two elements $\{1,-1\}$, its defining ideal in
$\K[\lambda_1,\lambda_2]$ being $\gid=(\lambda_2^2-1)$. 
The neutral element is $(0,1)$ while 
$(\bm_1,\bm_2)\cdot (\bl_1,\bl_2)^{-1}=(\bm_1-\bl_1,\bm_2\bl_2).$
We consider its action on $\K^2$ as translation parallel to the first
coordinate axis and reflection w.r.t. this axis. It is 
 defined  by the following polynomials of $\K[\lambda_1,\lambda_2,z_1,z_2]$:
$ g_1 = z_1+ \lambda_1 , \quad g_2=\lambda_2 z_2.$
\end{exy}

\begin{exy} \label{rotation:def} \textsc{rotation.}
Consider the special orthogonal group given by
$\gid= (\lambda_1^2+\lambda_2^2-1)\subset \K[\lambda_1, \lambda_2]$
with $e=(1,0)$ and $(\bm_1,\bm_2)\cdot (\bl_1,\bl_2)^{-1}=
 (\bm_1\bl_1+\bm_2\bl_2,\bm_2\bl_1-\bm_1\bl_2).$
Its linear action  on $\K^2$ is given by the following polynomials of
$\K[\lambda_1,\lambda_2,z_1,z_2]$:
$$ g_1 = \lambda_1 z_1-\lambda_2 z_2, \quad g_2=\lambda_2
z_1+\lambda_1z_2.$$
An element of the group acts as a rotation around the origin.
\end{exy}

\subsection{Graph of the action and orbits}\label{graph::orbits}

The \emph{graph of the action} is the image  
$\ova\subset\zva\times\zva$  of the map 
  $(\bl,\bz)\mapsto(\bz,g(\bl,\bz))$
that is defined on a dense open set of $\gva\times\zva$. We have
$\ova = \{(\bz,\bz') \;|\;  \exists \bl \in \gva \, \mbox{s.t.} \, \bz'=g(\bl, \bz)\}
\subset \zva\times\zva$.

We introduce a new set
of variables $Z=(Z_1, \ldots, Z_{\dimz})$ and the ideal
\( J=\gid+ (Z-g(\lambda,z)) \subset h^{-1}\K[\lambda,z,Z]\),
where $(Z-g(\lambda,z))$ stands for
$\left(Z_1-g_1(\lambda, z), \ldots\right.$, $\left. Z_{\dimz}-g_{\dimz}(\lambda,
z)\right)$.
 The set
$\ova$ is dense in its closure $\overline{\ova}$, and
$\overline{\ova}$ is the algebraic variety of the ideal:
\[ \oid = J\cap \K[z,Z]
=\left(\gid + (\,Z-g(\lambda, z)\,) \,\right)\cap \K[z,Z] .\]

Since $G$ is radical and unmixed dimensional so is $J$ because of 
the linearity in $Z$.
If $\gid=\bigcap_{i=0}^{\kappa} \gid^{(i)}$ is the prime decomposition
of $\gid$ then we have the following prime decomposition of $J$:
$$\left(\gid + (\,Z-g(\lambda, z)\,) \,\right)= \bigcap_{i=0}^{\kappa} \left(\gid^{(i)} +
  (\,Z-g(\lambda, z)\,) \,\right). $$ 
The prime ideal 
$\oid^{(i)} =  \left(\gid^{(i)} +  (\,Z-g(\lambda, z)\,) \,\right) \cap \K[z,Z]$
is therefore a component of $O$. 
The ideals  $\oid^{(i)}$, however, need not be all distinct. 

The set $\ova$ is symmetric: if
$(\bz,\bz')\in \ova$ then $(\bz',\bz)\in\ova$.
By the NullStellensatz the ideal $\oid$  is
also symmetric:  $p(Z,z)\in \oid$ if $p(z,Z)\in \oid$.
Since $J\cap\K[z]=(0)$, $O\cap\K[z]=(0)$ 
and therefore $\oid\cap \K[Z]=(0)$ also.

A set of generators, and more precisely a  Gr\"obner basis \cite{weisp},
  for $\oid \subset \K[z,Z]$ can be computed. 

\begin{propy} \label{gbforo}
Let $g'$ be the $\dimz$-tuple of numerators of $g$, that is 
$g'=hg=(h g_1, \ldots, h g_n) \in \left(\K[\lambda,z]\right)^n$. 
Consider a term order s.t.  $z \cup Z \ll \lambda \cup \{y\}$ where
$y$ is a new indeterminate.
If $Q$ is a Gr\"{o}bner basis for $\gid + (h\,Z - g') + (y h -1)$ according to
this term order then $Q \cap \K[z,Z]$ is a Gr\"{o}bner basis of $O$ according the
induced term order on    $z \cup Z$.
\end{propy}

\begin{proof}
Take $J' = (\gid + (Z-g))\cap \K[\lambda,z,Z]$ and note that
$J' = (\gid + (h\,Z - g'))\q h^\infty$  where $g'$ is the
numerator of $g$.
Given a basis  $\Lambda$ of $\gid$ and $g$
explicitly, a Gr\"obner basis of $\rhid$ is obtained thanks to
\cite[Proposition 6.37, Algorithm 6.6]{weisp}.
We recognize that $O$ is an elimination ideal of $J'$, namely
$O=J'\cap \K[z,Z]$. A Gr\"oner basis for $O$  is thus obtained by
\cite[Proposition 6.15, Algorithm 6.1]{weisp}.
\end{proof}

We mainly use the extension $\oid^e$ of $O$ in $\K(z)[Z]$.
If $Q$ is a Gr\"{o}bner basis of $O$ w.r.t. a term order $z\ll Z$ then
$Q$ is also a Gr\"{o}bner basis for $\oid^e$, for the term order
induced on $Z$ \cite[Lemma 8.93]{weisp}.
It is nonetheless often preferable to compute a Gr\"{o}bner basis of $\oid^e$
over $\K(z)$ directly (see \ex{derksen}).

The \emph{orbit} of $\bz\in\zva$ is the image $\ova_{\bz}$ of the rational map
$g_{\bz}\colon\gva\mapsto\zva$ defined by $g_{\bz}(\bl)=g(\bl,\bz)$. 
We then have the following specialization property (see for instance
\cite[Exercise 7]{cox}).

\begin{propy} \label{specialize}
Let $Q$ be a Gr\"{o}bner basis for $\oid^e$ for a given term order on
$Z$. There is a closed proper  subset $\mathcal{W}$ of
$\zva$ s.t. for $\bz\in \zva\setminus\mathcal{W}$ the image of $Q$
under the specialization $z \mapsto \bz$ is a Gr\"{o}bner basis 
 for the ideal whose variety is  the closure of the orbit of $\bz$.
\end{propy}

Therefore, for $\bz \in \zva\setminus\mathcal{W}$, the dimension of
the orbits of $\bz$  is equal to the dimension of
 $\oid^e \subset \K(z)[Z]$ \cite[Section 9.3, Theorem 8]{cox}.
In the rest of the paper  this dimension is denoted by $\dimo$.

\begin{exy} \label{scaling:graph} \textsc{Scaling}.
 Consider the group action of \ex{scaling:def}. 
The set of orbits consists of 1-dimensional punctured straight lines
through the origin and a single zero-dimensional orbit, the origin. 
By elimination on the ideal
$\rhid=( 1-\lambda_1\lambda_2, Z_1- \lambda_1 z_1, Z_2-\lambda_1 z_2)$
we  obtain $\oid=(z_1Z_2-z_2Z_1)$.
Take $\mathcal{W}$ to consist solely of the origin.  For
$\bz \in \zva\setminus \mathcal{W}$ the closure of the orbit of $\bz$
is  the algebraic variety of $(\bz_1Z_2-\bz_2Z_1)$
\end{exy}

\begin{exy} \label{translation:graph} \textsc{translation+reflection.}
Consider the group action of  \ex{translation:def}.
By elimination on the ideal
$\rhid = (\lambda_2^2-1, Z_1-z_1-\lambda_1, Z_2-\lambda_2 z_2)$ we
obtain $\oid=(Z_2^2-z_2^2)$. 
The orbit of a point $\bz=(\bz_1,\bz_2)$ with $\bz_2\neq 0$ 
consists of two lines parallel to the first
coordinate axis, while the latter
 is the orbit of all points with $\bz_2= 0$
\end{exy}

\begin{exy} \label{rotation:graph} \textsc{rotation.}
Consider the group action of \ex{rotation:def}. 
The orbits consist of the origin and the circles with the origin as
center. By elimination on the ideal
$\rhid = (\lambda_1^2+\lambda_2^2-1, Z_1-\lambda_1 z_1+\lambda_2 z_2,
Z_2-\lambda_2 z_1-\lambda_1z_2)$ we obtain $\oid=(Z_1^2+Z_2^2-z_1^2-z_2^2)$.
\end{exy}

\subsection{Rational invariants} \label{rational}

We construct  a finite set of  generators for the field of rational invariants.  Our construction 
brings out a simple algorithm to rewrite any rational invariant in
terms of them. The required operations are restricted to computing a Gr\"{o}bner basis and
normal forms. Those are implemented in most  computer algebra systems. We provide a comparison with related results in \cite{beth99,rosenlicht56,vinberg89}.

\begin{defy} A rational function $r\in \K(z)$ is a 
\emph{rational invariant} if \phantom{o} $r(g(\lambda,z)) = r(z) \mod G.$
\end{defy}

The set of rational invariants forms a 
field\footnote{Though we do not use this fact but rather retrieve it
  otherwise, it is worth noting that, as a subfield of $\K(z)$,
  the field of rational invariants  is always finitely generated \cite{waerden71}.}
$\ratif$. 
We show that the coefficients of the Gr\"obner basis for
$\oid^e$ are invariant and generate $\ratif$. The basis is computed
using \rprop{gbforo}.

\begin{lemy} \label{invariance} If $q(z,Z)$ belongs to $\oid$ then
$q(g({\bl},z),Z)$ belongs to $\oid^e$ for
all $\bl\in \gva$.
\end{lemy}

\begin{proof} A point $(\bz,\bz')\in \ova$ if there exists $\bm\in\gva$
s.t. $\bz'=g(\bm,\bz)$. Then for a generic $\bl\in \gva$, $\bz'=
g(\bm \cdot \bl^{-1}, g(\bl,\bz))$. Therefore $(g(\bl,\bz),\bz')\in\ova$.
Thus if $q(z,Z)\in \oid$ then $q(g({\bl},\bz),\bz')=0$ for all
$(\bar{z},\bz')$ in $\ova$. By the Hilbert NullStellensatz the
numerator of $q(g({\bl},z),Z)$  belongs to $\oid$ and therefore
$q(g({\bl},z),Z)\in \oid^e$.
\end{proof}

Following \cite[Definition 5.29]{weisp}  a set of polynomials is
reduced, for a given term order, if   the leading coefficients of the elements  are equal
to  $1$ and each element is in normal form with respect to the others.
 Given a term order on 
 $Z$, a polynomial ideal in $\K(z)[Z]$ has a unique reduced
Gr\"{o}bner basis \cite[Theorem 5.3]{weisp}.

\begin{theoy} \label{invgb} The reduced Gr\"{o}bner basis of $\oid^e$
with respect to any term order on $Z$ consists of polynomials in
$\ratif[Z]$.
\end{theoy}

\begin{proof}
Let $Q=\{q_1, \ldots, q_\kappa\}$ be the reduced
Gr\"{o}bner basis of $\oid^e$ for a given term order on $Z$.
By \rlem{invariance} $q_i(g(\bl,z), Z)$ belongs to $\oid^e$. It has
the same support\footnote{The support here is the set of terms in $Z$ with non zero
coefficients.} as $q_i$. As $q_i(g(\bl,z), Z)$ and $q_i(z, Z)$ have
the same leading monomial, $q_i(g(\bl,z), Z) - q_i(z,Z)$ is in
normal form with respect to $Q$. As this difference belongs to ${\oid}^e$, it
must be $0$.  The coefficients of $q_i$ are therefore invariant.
\end{proof}

Let us note the  construction
of a generating set of rational invariants proposed by Rosenlicht
\cite{rosenlicht56}. In the paragraph before
Theorem 2,  Rosenlicht points out that the coefficients of the Chow form of
$\oid^e$  over  $\K(z)$  form  a set of separating rational invariants.
By \cite[Theorem 2]{rosenlicht56} or \cite[Lemma 2.1]{vinberg89} this
set is generating for $\ratif$.

Vinberg and Popov showed the existence of a subset of $\ratif[Z]$  that
generates  $\oid^e$  \cite[Lemma 2.4]{vinberg89}.  We propose the
construction of such a set. They showed furthermore  
that the set of the  coefficients   of such a family of generators
\emph{separates generic orbits}
\cite[Theorem 2.3]{vinberg89} and therefore generates  $\ratif$
\cite[Theorem 2]{rosenlicht56},\cite[Lemma 2.1]{vinberg89}.
From those results we deduce that the set of coefficients of a
reduced Gr\"obner basis of $\oid^e$ generates $\ratif$. The next
theorem provides an alternative proof of this result, providing
additionally  a rewriting algorithm. To prove generation we
indeed  exhibit an algorithm that allows to rewrite any rational
invariant in terms of the coefficients of a reduced Gr\"obner
basis.

In the case of linear actions
 M{\"u}ller-Quade and Beth \cite{beth99} showed
that the coefficient of the Gr\"obner basis of $\oid^e$ generate
the field of rational invariants. Their proof is based  on 
more general results about the characterization of 
subfields of $\K(z)$ obtained in \cite{quade99}. 
Our approach is  quite different and more direct. 
The rewriting algorithm we propose,  although it was obtained
independently, is  nonetheless reminiscent of \cite[Algorithm 1.10]{quade99}.

\begin{lemy} \label{lolita} Let $\frac{p}{q}$ be a rational invariant,
  $p,q\in \K[z]$.   Then \( p(Z)\,q(z)- q(Z)\, p(z) \in \oid .\)
\end{lemy}

\begin{proof}
Since $\frac{p}{q}$ is an invariant
$\frac{p(\bz)}{q(\bz)}=\frac{p( g(\bl,\bz))}{q(g(\bl,\bz))}$
for all $(\bl,\bz)$ where this expression is defined.
Thus  \(a(\bz',\bz)=p(\bz')\,q(\bz)- q(\bz')\, p(\bz)=0\)
for all $(\bz,\bz')$ in
$\ova=\left\{(\bz,{\bz}')\,|\,
\exists\bl\in \gva \mbox{ s.~t. }{\bz}'=g(\bl,\bz)\right\}\subset\zva\times\zva$.
In other words the polynomial
\(a(Z,z)= p(Z)\,q(z)- q(Z)\, p(z)\in \K[Z,z] \)
is zero at each point of $\ova$.
Since the algebraic variety of $\oid$ is the closure $\bar\ova$ of
$\ova$ and that $\ova$ is dense in $\bar{\ova}$
we can conclude that  \(a(Z,z)\in \oid\) by  Hilbert Nullstellensatz.
\end{proof}

Assume a polynomial ring over a field is endowed with a given  term order.
A  polynomial $p$ is in \emph{normal form} w.r.t.  a set $Q$ of
polynomials if   $p$ involves no term that is a multiple of a
leading term of an   element in $Q$.
A \emph{reduction} w.r.t. $Q$  is an algorithm
that given $p$ returns a polynomial $p'$ in normal form w.r.t. $Q$
s.t. $p = p' + \sum_{q\in Q} a_q \, q$ and no leading term of
any $a_q \, q$ is larger than the leading term of $p$.
Such an algorithm is detailed in \cite[Algorithm 5.1]{weisp}. It
consists in rewriting the terms that are multiple of the leading terms
of the elements of $Q$ by polynomials involving only terms that are lower.
Note that if the leading coefficients of $Q$ are $1$ then no division
occurs.
When $Q$ is a Gr\"{o}bner basis w.r.t. the given term order,
the reduction of a polynomial $p$ is unique in the sense that $p'$ is
then the only polynomial in normal form w.r.t. $Q$
in the equivalence class $p+(Q)$.

\begin{theoy} \label{rewrite}
Consider $\{r_1,\ldots, r_\kappa\} \in \ratif$ the coefficients
of a reduced Gr\"{o}bner  basis $Q$ of $\oid^e$.
Then $\ratif=\K(r_1,\ldots, r_\kappa)$ and we can
rewrite any rational invariant
$\frac{p}{q}$, with $p,q\in \K[z]$,
in terms of those as follows.

Take a new set of indeterminates $y_1, \ldots, y_\kappa$ and
consider the set $Q_y \subset \K[y,Z]$ obtained from $Q$
 by substituting $r_i$ by  $y_i$.
Let $a(y, Z)=\sum_{\alpha\in\Ni^{\dimz}}a_{\alpha}(y)Z^\alpha$ and
$b(y,Z)=\sum_{\alpha\in\Ni^n}a_{\alpha}(y)Z^\alpha$ in $\K[y,Z]$
be the reductions\footnote{For those reductions in $\K[y,Z]$ the term
  order on $Z$ is extended to a block order $y \ll Z$ so that the set
  of leading term of $Q_y$ is equal to the set of leading terms of $Q$.}
of $p(Z)$ and $q(Z)$ w.r.t. $Q_y$.
There exists $\alpha\in \Ni^{\dimz}$ s.t. $b_{\alpha}(r)\neq 0$ and
for any such
$\alpha$ we have $\frac{p(z)}{q(z)}=\frac{a_\alpha(r)}{b_{\alpha}(r)}$.
\end{theoy}

\begin{proof} It is sufficient to prove the second part of the
  statement.
The Gr\"{o}bner basis $Q$ is reduced and therefore
  monic. The sets of leading monomials of $Q$ and of $Q_y$ are equal.
If $a(y,Z)$ is the reduction of $p(Z)$ w.r.t. $Q_y$ then
$a(r,Z)$, obtained by substituting back $y_i$ by $r_i$,
is the normal form of $p(Z)$ w.r.t. $Q$.
Similarly for  $b(y,Z)$ and  $q(Z)$.

As $\oid^e\cap\K[Z]=(0)$, neither $p(Z)$ nor $q(Z)$ belong
to $\oid^e$ and therefore both $a(r,Z)$ and $b(r,Z)$ are different
from $0$. By \rlem{lolita}  $q(z)p(Z) \equiv p(z)q(Z) \mod \oid^e$
and thus the normal forms of the two polynomials modulo $\oid^e$ are equal:
 $q(z)\,a(r,Z)=p(z)\,b(r,Z)$. 
Thus $a(r,Z)$ and $b(r,Z)$ have the same support and
this latter is non empty since $a,b\neq 0$.
 For each $\alpha$ in this common support, we have
$q(z) a_{\alpha}(r) = p(z) b_{\alpha}(r)$
and therefore $\frac{p(z)}{q(z)} = \frac{a_{\alpha}(r)}{b_{\alpha}(r)}$.
\end{proof}

\begin{exy} \label{scaling:rational} \textsc{Scaling}.
We consider the group action  given in \ex{scaling:def}.
A reduced Gr\"{o}bner basis of $\oid^e$
is $Q=\{ Z_2 -\frac{z_2}{z_1}Z_1\}$.  By \theo{invgb},
$\K(z_1,z_2)^G=\K(\frac{z_2}{z_1})$.

Let $p = z_1^2+4z_1z_2+z_2^2$ and $q=z_1^2-3z_2^2$.
We can check that $p/q$ is a
rational invariant and we set up to write $p/q$ as a rational function
of $r=z_2/z_1$. To this purpose
consider $P=Z_1^2+4Z_1Z_2+Z_2^2$ and $Q=Z_1^2-3Z_2^2$ and compute  their
normal forms  $a$ and $b$ w.r.t.  $\{Z_2 - y\,Z_1\}$
 according to a term order where $Z_1<Z_2$.
We have  $a=(1+4y+y^2)Z_1^2$ and $b=(1-3y^2)Z_1^2$. Thus
\[  \frac{z_1^2+4z_1z_2+z_2^2}{z_1^2-3z_2^2} = \frac{1+4r+r^2}{1-3r^2}
\hbox{ where } r= \frac{z_2}{z_1} \]
\end{exy}

\begin{exy} \label{translation:rational} \textsc{translation+reflection}.
We consider the group action  given in \ex{translation:def}.
A reduced Gr\"{o}bner basis of $\oid^e$
is $Q=\{ Z_2^2 -z_2^2\}$. By \theo{invgb},
$\K(z_1,z_2)^G=\K(z_2^2)$.
\end{exy}

\begin{exy} \label{rotation:rational} \textsc{Rotation}.
We consider the group action  given in \ex{rotation:def}.
 A reduced Gr\"{o}bner basis of $\oid^e$
 is $Q=\{ Z_1^2+Z_2^2 -(z_1^2+z_2^2)\}$. By \theo{invgb},
$\K(z_1,z_2)^G=\K(z_1^2+z_2^2)$.
\end{exy}

\section{Cross-section and rational invariants} 
\label{crossection}

Given a cross-section we  construct a  generating set of
rational invariants endowed with a rewriting algorithm. The method is
the same as the one presented in previous section but applies to only a
section of the graph. In  previous section we considered an ideal of
the  dimension of the generic orbits.Here we  consider  a zero dimensional ideal.
This is computationally advantageous when Gr\"obner bases are needed.

We use  Noether normalization to prove the existence of a
 cross-section.
The construction thus relies on selecting elements of in an open
subset of a certain affine space. 
Therefore the  construction   does not entail a deterministic algorithm for
the computation of rational invariants. 
Yet the freedom of choice is extremely fruitful in applicative
examples.

Though the presentation is done with an algebraically closed
field $\K$ that is therefore infinite, the construction 
is  meant to be realized in characteristic zero 
(i.e. over $\Qi$) or over a sufficiently large field. 

\subsection{Cross-section} \label{acs}

Geometrically speaking a \emph{cross-section  of degree $d$} is a
variety that intersects generic orbits in $d$ simple points.
We give a definition in terms of ideals for it is closer to the actual
computations. We give its geometric content in a proposition afterward.

\begin{defy}\label{dcs} Let $\sid$ be a prime ideal of $\K[Z]$  of
  complementary dimension to the generic orbits, i.e. if $O^e$ is of
  dimension ${\dimo}$ then  $\sid$ is of codimension $\dimo$.
  $\sid$  defines a
  \emph{cross-section} to the orbits of the rational action
  $g:\gva\times\zva\rightarrow \zva$ if the ideal
  $\mfid^e = \oid^e + \sid$ of $\K(z)[Z]$ is radical and zero dimensional.
  If $d$ is the dimension of  $\K(z)[Z]/\mfid^e$ as a $\K(z)$-vector
space, we say that $\sid$ defines a \emph{cross-section of degree $d$}.
\end{defy}

Indeed
the algebra  $\K(z)[Z]/\mfid^e$ is a finite $\K(z)$-vector
space since  $I^e$ is zero dimensional \cite[Theorem 6.54]{weisp}.
A basis for it is provided by the terms in $Z$ that are not
multiple of the leading terms of a Gr\"{o}bner basis of $\mfid^e$
\cite[Proposition 6.52]{weisp}.
Let us note here that  an ideal of $\K(z)[Z]$ is zero dimensional iff
any Gr\"obner basis of it has an element whose leading term is
$Z_i^{d_i}$, for  all $1\leq i\leq n$ \cite[Theorem 6.54]{weisp}.

The cross-section is thus the variety $\sva$ of $\sid$.
The geometric properties of this variety are explained by the following proposition.

\begin{propy} \label{transverse}
Let $\sid$ define a cross-section $\sva$ of degree $d$. 
There is a closed set $\mathcal{S} \subset\zva$
s.t. the closure of the orbit of any $\bz\in\zva\setminus\mathcal{S}$
intersects $\sva$ in $d$  simple points.
\end{propy}

\begin{proof}
Let $Q$ be a reduced Gr\"{o}bner basis for $\mfid^e=\oid^e+\sid$.
Similarly to \rprop{specialize}, the image $Q_{\bz}$ of $Q$ under the
specialization $z\mapsto\bz$ is a Gr\"{o}bner basis for
$\oid_{\bz}+\sid$ in $\K[Z]$
for all $\bz$ in $\zva$ outside of a closed set $\mathcal{W}$.
Thus $I_{\bz}=\oid_{\bz}+\sid$ is zero dimensional and the dimension of
$\K[Z]/I_{\bz}$ as a vector space over $\K$ is $d$.

By the Jacobian criterion for regularity and the prime avoidance
theorem \cite[Corollary 16.20 and Lemma3.3]{eisenbud} there is a
$n\times n$ minor $f$ of the Jacobian matrix of $Q$ that is
not included in any prime divisor of $\mfid^e$. Therefore $f$ is not a
zero divisor in $\K(z)[Z]/\mfid^e$ which is a product of
fields. There exists  thus $f'\in\K(z)[Z]$
s.t. $f\,f'\equiv 1 \mod \mfid^e$.

Provided that $\bz$ is furthermore chosen so that the
denominators of $f$ and $f'$ do not vanish, $f$  specializes into
a $n\times n$ minor ${f}_{\bz}$ of the Jacobian
matrix of $Q_{\bz}$ and we have
$f_{\bz}\,f_{\bz}'\equiv 1 \mod \mfid_{\bz}$
for the specialization $f_{\bz}'$ of $f'$.
So $f_{\bz}$ belongs to no prime divisors of  $\mfid_{\bz}$
and thus $\mfid_{\bz}$ is radical \cite[Corollary 16.20]{eisenbud}.
We take  $\mathcal{S}$ to be the union of $\mathcal{W}$ with
the algebraic set associated to the product of the denominators of $f$
and $f'$. That the number of points of intersection is $d$ is shown by
\cite[Proposition 2.15]{eisenbud}.
\end{proof}

That property shows that the cross-sections of degree $d=1$ and $d>1$
are  respectively the sections and the quasi-sections 
defined in \cite[Paragraph 2.5]{vinberg89}.
The existence of
quasi-section  is insured by \cite[Proposition 2.7]{vinberg89}, while
a criterion for the existence of a section is described in
\cite[Paragraph 2.5 and 2.6]{vinberg89} 
Our terminology elaborates on the one used in 
\cite{rosenlicht56} and \cite{olver99}. 

The discussion of \cite[Section 2.5]{vinberg89}
shows  that $\K(\sva)$ is isomorphic to $\ratif$ when $\sva$ 
is a cross-section of degree 1.
If $\sva$ is a cross-section of degree $d>1$ then
$\K(\sva)$ is an algebraic extension of $\ratif$ of degree $d$.
In \sec{invariantization} we shall come back to those points
with a constructive angle that relies on the choice of a cross-section.
The viewpoint  adopted here is indeed the geometric intuition  of the moving
frame construction in \cite{olver99}: almost any algebraic variety of
complementary dimension provides a cross-section (of some degree).

The existence of a cross-section  is proved  by Noether
normalization theorem and is linked to an alternative
definition of the dimension of an ideal \cite[Section 6.2]{shafarevich}.

\begin{theoy} \label{cross:existence}
A  linear cross-section to the orbit is associated to each point of an
 open set of $\K^{\dimo(n+1)}$, where $\dimo$ is the dimension of the generic
 orbits and $\dimz$ the dimension of $\zva$.
\end{theoy}
\vskip2mm
\begin{proof}
Assume that a Gr\"obner basis $Q$ of $\oid^e$ w.r.t. a term order
$Z_1,\ldots, Z_{\dimo} \ll Z_{\dimo+1},\ldots, Z_{\dimz}$ is s.t.
 an element  of $Q$
has leading term $Z_i^{d_i}$ for some $d_i\in \Ni\setminus\{0\}$
 for all $\dimo+1\leq i\leq \dimz$
and  there is no element of $Q$ independent of $\{ Z_{\dimo+1},\ldots, Z_{\dimz}\}$.
Then $Q$ is a Gr\"obner basis for the extension of $\oid^e$ to
$\K(z)(Z_1,\ldots,Z_{\dimo})[Z_{\dimo+1}, \ldots, Z_{\dimz}]$
\cite[Lemma 8.93]{weisp}.
For $(a_1,\ldots, a_{\dimo})$ in an open set of  $\K^\dimo$
the specialization $Q_a \subset \K[Z_{\dimo+1}, \ldots, Z_{\dimz}]$
of $Q$ under $Z_i\mapsto a_i$ is a Gr\"obner
basis \cite[Exercise 7]{cox}.
Therefore $Q_a \cup \{Z_1-a_1, \ldots, Z_{\dimo}-a_{\dimo}\}$ is a
Gr\"obner basis by Buchberger's criteria \cite[Theorem 5.48 and 5.66]{weisp}.
It is a Gr\"obner basis of a zero dimensional ideal \cite[Theorem 6.54]{weisp}.
We can thus take $\sid$ to be generated by $\{Z_1-a_1, \ldots, Z_{\dimo}-a_{\dimo}\}$.

We can always retrieve the situation assumed above by a
change of variables thanks to  Noether normalization
theorem \cite[Theorem 3.4.1]{greuel}.
Inspecting the proof we observe that we can choose a change of variables
given by a matrix $(m_{ij})_{1\leq i,j\leq \dimz}$
with the vector of entries  $m_{ij}$ in $\K^{{\dimz}^2}$ outside of some
algebraically closed set.
The set $\{ a_i - \sum_{1\leq j \leq \dimz} m_{ij} Z_j\;|\; 1\leq i\leq \dimo\}$
thus defines a cross-section.~\end{proof}

The choice of a cross section introduces a non deterministic aspect 
to the algebraic construction proposed in next section. 
An analysis of the probability of success in characteristic $0$ 
would be based on the measure of a correct test sequence
\cite[Theorem 3.5 and 3.7.2]{giusti93a},
\cite[Section 3.2]{giusti93b}, \cite[Section 4.1]{krick01}. 

We can computationally test if $\sva$ is a cross-section by checking
the properties of $\mfid^e= \left(G+ P+ (Z-g(\lambda,z))\right)\cap\K(z)[Z]$, 
starting with the computation of its Gr\"obner basis. 
It is nonetheless worth performing the preliminary necessary test 
of transversality detailed in \sec{smooth::section}. It relies on
computing the rank of a matrix.

\begin{propy} \label{mfidcomp}
Assume that $\sid\subset \K[Z]$ defines a
  cross-section and that $\oid = \bigcap_{i=0}^{\tau} \oid^{(i)}$ is
  the prime decomposition of $O$. Then
 \[ \oid + \sid = \bigcap_{i=0}^{\tau}(\oid^{(i)}+\sid) \quad \hbox{
   and } \quad
 (\oid^{(i)}+\sid) \cap \K[Z] = \sid.\]
\end{propy}

\begin{proof}
We can easily check that
$ \bigcap_{i=0}^{\tau}(\oid^{(i)}+\sid) \subset \oid + \sid $ because
$\oid+\sid$ is radical. The converse inclusion is trivial.

For the second equality, note first that $\sid \subset
(\oid^{(i)}+\sid)\cap\K[z,Z]$. The projection of the variety of
$\oid^{(i)} \subset \zva\times\zva$ is thus contained in $\sva$.
We show that the projection is exactly $\sva$.
We can assume that the numbering is such that
$\oid^{(i)}=  (\left(\gid^{(i)} +  (\,z-g(\lambda,  Z)\,) \,\right)
\cap \K[z,Z]$ where $\gid^{(i)}$ is a minimal prime of $\gid$ (see
\sec{groupaction}).
By \hyp{groupaction:hyp}, 
for any $\bz$ in $\zva$ and therefore in $\sva$,
there exists $\bl$ in the variety of $\gid^{(i)}$
s.t. $g(\bl, \bz)$ is defined.
Above each point of $\sva$ there is a point in the variety of
$\oid^{(i)}$.
\end{proof}

\subsection{Rational invariants revisited} \label{mango}

The following theorems  provide a construction of a generating set of
rational invariants together  with an algorithm to  rewrite any rational invariant in terms of
generators. The method is the same as in \sec{rational} but applied to
the ideal $I^e$ rather than to $\oid^e$. The computational
advantage comes from the fact that  $\mfid^e$
is zero dimensional.

If $G$ is a prime ideal we can actually choose a coordinate
cross-section that is $\sid$ can be taken as the ideal generated by a set of
the following form:
$\{Z_{j_1}-\alpha_1, \ldots, Z_{j_{\dimo}}-\alpha_{\dimo}\}$ for
$(\alpha_1, \ldots, \alpha_{\dimo})$ in $\K^{\dimo}$. In this case
we can  remove $r$ variables for the computation.

\begin{theoy}\label{cocorico}
 The reduced Gr\"{o}bner basis of $\mfid^e$
with respect to any term ordering on $Z$ consists of polynomials in
$\ratif[Z]$.
\end{theoy}

\begin{proof}
The union of a reduced Gr\"{o}bner basis  of
$\oid^e$ and  $\sid$ forms a generating set for
$\mfid^e=\oid^e + \sid$. The coefficients of a  basis for   $\sid$
are  in $\K$, while the coefficients of a reduced basis for $\oid^e$
belong to  $\ratif$
due to \theo{invgb}.  Since the coefficients of a generating set for
$\mfid^e$ belong to  $\ratif$, so do the  coefficients of the reduced
Gr\"{o}bner basis with respect to any term ordering.
\end{proof}

\begin{theoy} \label{rewrite2}
 Consider  \phantom{o}  $\{r_1,\ldots, r_\kappa\} \in \ratif$\phantom{o}  the coefficients 
of a reduced  Gr\"{o}bner  basis $Q$ of $\mfid^e$.
Then $\ratif=\K(r_1,\ldots, r_\kappa)$ and we can
rewrite any rational invariant
$\frac{p}{q}$, with $p,q\in \K[z]$ relatively prime,
in terms of those as follows.

Take a new set of indeterminates $y_1, \ldots, y_\kappa$ and
consider the set $Q_y \subset \K[y,Z]$ obtained from $Q$
 by substituting $r_i$ by  $y_i$.
Let $a(y, Z)=\sum_{\alpha\in\Ni^n}a_{\alpha}(y)Z^\alpha$ and
$b(y,Z)=\sum_{\alpha\in\Ni^n}a_{\alpha}(y)Z^\alpha$ in $\K[y,Z]$
be the  reductions of $p(Z)$ and $q(Z)$ w.r.t. $Q_y$.
There exists $\alpha\in \Ni^m$ s.t. $b_{\alpha}(r)\neq 0$ and for any such
$\alpha$ we have
$\frac{p(z)}{q(z)}=\frac{a_\alpha(r)}{b_{\alpha}(r)}$.
\end{theoy}

\begin{proof}
We can  proceed just as in the proof of \theo{rewrite}.
We only need to argue additionally that if
$r=\frac{p}{q}\in\ratif$, $p$ and $q$ being relatively prime, then
$p(Z),q(Z)\notin\mfid^e$. We prove the result for $p$, the case of $q$
being similar.

By hypothesis  $p(z)\, q(g(\lambda,z)) \equiv q(z)\, p(g(\lambda,z)) \mod \gid$.
Since $p$ and $q$ are relatively prime, $p(z)$ divides $p(g(\lambda,z))$
modulo $\gid$, that is there exists $\alpha \in h^{-1}\K[z,\lambda]$
s.t.  $p(g(\lambda,z))\equiv \alpha(\lambda,z) \, p(z)\mod \gid$.
 Therefore if $p$ vanishes at
$\bz\in\zva$, then it vanishes on $\ova_{\bz}$.
Thus if $p\in\sid$, or equivalently if $p$ vanishes on $\sva$, 
it vanishes on an open subset of $\zva$ (\rprop{transverse}). 
So $p$ must be zero. This is not the case and thus
$p\notin\sid$. Since $\mfid^e\cap \K[Z]=\sid$, it is the case that
$p(Z)\notin I^e$
\end{proof}

When $\sid$ defines a cross-section
of degree $1$, the rewriting trivializes into a  \emph{replacement}.
Indeed, if  the dimension of $\K(z)[Z]/\mfid^e$ as a $\K(z)$ vector
space is $1$
then, independently of the chosen term order, the reduced Gr\"{o}bner
basis $Q$ for $I^e$ is given by $\{ Z_i -r_i(z) \,|\, 1\leq i\leq n\}$ where
the $r_i\in \ratif$. In view of \theo{rewrite2}
 $\ratif=\K(r_1, \ldots, r_\dimz)$ and any
rational invariant $r(z)\in \ratif$ can be rewritten in terms of $r_i$
by replacing $z_i$ by $r_i$:
$$r(z_1,\ldots,z_\dimz)= r(\, r_1(z), \ldots, r_{\dimz}(z)\,),
\quad \forall r\in \ratif.$$
In the next section we generalize this replacement to cross-section of
any degree by introducing some special \emph{algebraic invariants}.

\begin{exy}\label{scaling:cross} \textsc{scaling}.
We carry on with the action considered in
\ex{scaling:def} and \ref{scaling:rational}.

Choose  $\sid=(Z_1-1)$.  A reduced Gr\"{o}bner basis of
$\mfid^e$ is given by \(\{Z_1-1, Z_2 -\frac{z_2}{z_1}\}\).
We can see that \theo{cocorico} is verified and that $\sid$ defines a
cross-section of degree 1.
By \theo{rewrite2} we know that
$r=z_2/z_1$ generates the field of rational invariants $\ratif$. 
In this situation, the cross section is of degree $1$ and we see that
the rewriting algorithm of \theo{rewrite2}
is a simple replacement.
For all $p\in \ratif$ we have $p(z_1,z_2)=p(1,r)$.
\end{exy}

\begin{exy}\label{translation:cross} \textsc{translation+reflection}.
We carry on with the action considered in
\ex{translation:def} and \ref{translation:rational}.

Choose  $\sid=(Z_1-Z_2)$ to define the cross-section. A reduced Gr\"{o}bner basis of
$\mfid^e$ is given by \(\{Z_1-Z_2, Z_2^2 -z_2^2\}\).
The cross-section is thus of degree 2. 
\end{exy}

\begin{exy}\label{rotation:cross} \textsc{rotation}.
We carry on with the action considered in
\ex{rotation:def} and \ref{rotation:rational}.

Choose  $\sid=(Z_2)$.  The reduced Gr\"{o}bner basis of
$\mfid^e$ w.r.t. any term order is  \(\{Z_2, Z_1^2-(z_1^2+z_2^2)\}\).
We can see that \theo{invgb} is verified and that $\sid$ defines a
cross-section of degree 2. By \theo{rewrite2} we know that
$r=z_1^2+z_2^2$ generates the field of rational invariants 
$\ratif$. 
In this situation, the rewriting algorithm of \theo{rewrite2}
consists in substituting $z_2$ by $0$ and $z_1^2$ by $r$.
\end{exy}

\section{Replacement invariants and invariantization} 
   \label{invariantization}

Given a cross-section $\sva$ of degree $d$
we introduce $d$ distinct $\dimz$-tuples  
of  elements that are algebraic over the field of  rational invariants.
Each $n$-tuple has an important \emph{replacement} property: any
rational invariant can be rewritten in terms of its components by a simple
substitution of the variables by the corresponding elements  from the tuple.

The replacement invariants are used  to define a process of \emph{invariantization}, that is  
a projection  from the algebraic  functions  onto the field
of algebraic invariants.  This projection can be explicitly computed
by algebraic  elimination. It gives a constructive approach to the
isomorphism $\overline{\K(\sva)}\cong\algif$. 

\subsection{Replacement invariants}\label{repl}

Let $\sva$ be a cross-section of degree $d$ 
defined by a prime ideal  $\sid$ of $\K[Z]$. 
The field of rational  functions on $\sva$ is denoted by
$\K(\sva)$. It is the fraction field of the integral 
domain $\K[Z]/\sid=\K[\sva]$. 
We introduce $d$ replacement invariants associated to $\sva$. We use them to show that $\K(\sva)$ is an algebraic
extension of degree $d$ of the field of rational invariants
$\ratif$.

\begin{defy} An \emph{algebraic invariant} is an element of the
algebraic closure $\algif$ of $\ratif$.
\end{defy}

A reduced Gr\"obner basis   $Q$ of  $\mfid^e=\oid^e+\sid$
is contained in $\ratif[Z]$ (\theo{cocorico}) and therefore is a 
reduced Gr\"obner basis  of $\mfid^G=\mfid^e\cap\ratif[Z]$.
The dimension of $\ratif[Z]/\mfid^G$ as a $\ratif$-vector space 
is therefore equal to the dimension $d$ of $\K(z)[Z]/\mfid^e$ as a
$\K(z)$-vector space. 
Consequently the ideal $\mfid^G$ has  $d$  zeros 
$\ninv=(\ninv_1, \ldots, \ninv_n)$ with
$\xi_i\in \algif$ \cite[Proposition 2.15]{eisenbud}.
We call such a tuple  $(\xi_1, \ldots, \xi_n)$ a  $\algif$-zero of
$\mfid^G$. 
A $\algif$-zero of $\mfid^G$ is a  $\algif$-zero of
$\mfid^e$ and conversely.

\begin{defy} 
  A replacement invariant is a
 $\algif$-zero of $\mfid^G=\mfid^e\cap\ratif[Z]$, i.e. a $\dimz$-tuple 
$\ninv=(\ninv_1,\ldots, \ninv_n)$ of algebraic invariants that forms a
zero of  $\mfid^e$.
\end{defy}

 Thus  $d$ replacement
invariants $\ninv^{(1)}, \ldots,\ninv^{(d)}$ are associated to a cross-section of degree $d$.
The name owes to next theorem which  can be compared with Thomas replacement 
theorem discussed in \cite[page 38]{olver99}.

\begin{theoy} \label{replacement}
Let $\ninv=(\ninv_1,\ldots,\ninv_\dimz)$ be a replacement invariant.
If $r \in \ratif$ then \(r(z_1,\ldots, z_\dimz)=r(\ninv_1,\ldots, \ninv_\dimz)\) in $\algif $.
\end{theoy}

\begin{proof} Write $r = \frac{p}{q}$ with $p,q$ relatively prime. 
By \rlem{lolita}, $p(z)\,q(Z)- q(z)\,p(Z) \in \oid^e \subset \mfid^e$ and
therefore $p(Z)- \frac{p(z)}{q(z)}\, q(Z)=p(Z)-r(z) \, q(Z)
\in\mfid^e$. 
Since $\ninv$ is a zero of $\mfid^e$, we have 
$p(\ninv)-r(z) \, q(\ninv)=0$.  In the proof of
\theo{cocorico} we saw that $p(Z),q(Z)$ can not belong to $P$ and
therefore cannot be zero divisors modulo $\mfid^e$. Thus 
$q(\ninv)\neq 0$ and the  conclusion follows.
\end{proof}

The field $\K(\ninv)$, 
for any replacement invariant $\ninv$, 
is an algebraic extension of $\ratif$. Indeed
 $\ratif\subset\K(\ninv)$ and $\ninv$ is algebraic over
$\ratif$. This leads to the following results.

\begin{lemy} \label{IGprime}$I^G=\mfid^e\cap \ratif[Z]$ is a prime ideal of $\ratif[Z]$.  
\end{lemy}

\begin{proof}
Let $I^{(1)}$ and $I^{(2)}$ be prime divisors of $\mfid^G$ in $\ratif[Z]$
and consider  replacement invariants $\ninv^{(1)}$ and $\ninv^{(2)}$ that are
 $\algif$-zeros of $I^{(1)}$ and $I^{(2)}$ respectively. Due to
 \theo{replacement} $\K(\ninv^{(i)})=\ratif(\ninv^{(i)})$. There is 
 therefore   a $\ratif$-isomorphism  $\ratif[Z]/I^{(i)}\cong\K(\ninv^{(i)})$ 
 for $i=1$ or $2$. On the other hand we have  $\K(\ninv^{(i)})\cong\K(\sva)$
since $P$ is the ideal of all relationships on the components of $\ninv^{(i)}$ 
over $\K$ (\rprop{mfidcomp}).
 Thus
\[ \ratif[Z]/I^{(1)}\cong\K(\ninv^{(1)})\cong \K(\sva)\cong\K(\ninv^{(2)})\cong
\ratif[Z]/I^{(2)}. \]
We have an isomorphism between 
$\ratif[Z]/I^{(1)}$ and $\ratif[Z]/I^{(2)}$ that leaves  $\ratif$
fixed and
maps the class of $Z$ modulo $ I^{(1)}$ to the class of $Z$ modulo 
$I^{(2)}$. Therefore $I^{(1)}=I^{(2)}$ so that $\mfid^G$ is prime.
\end{proof}

\begin{theoy} \label{d-ext}  The field   $\K(\sva)$ is an algebraic
  extension of $\ratif$ of degree~$d$,  the degree of the
  cross-section  $\sva$. 
  \end{theoy}

  \begin{proof}
For any replacement invariant $\ninv$ we have 
\(\ratif[Z]/\mfid^G\cong\K(\ninv)\cong \K(\sva)\). 
Since the dimension of $\ratif[Z]/\mfid^G$ as $\ratif$-vector space  is $d$, 
the field $\K(\sva)$ is an algebraic extension of $\ratif$ of degree~$d$.
\end{proof}

In particular if $\sva$ is a cross-section of degree one we have 
$\K(\sva)\cong \ratif$. In all cases we have the isomorphism
$\overline{\K(\sva)}\cong\algif$ obtained in  \cite[Section~2.5]{vinberg89} by
different means.

\begin{exy}\label{scaling:algebraic} \textsc{scaling}.
Consider the multiplicative group from 
\ex{scaling:def}, \ref{scaling:graph}, \ref{scaling:rational}. 
We considered the cross-section of degree 1  
defined by $\sid=(Z_1-1)$.  There is single replacement invariant 
$\ninv=(1,\frac{z_2}{z_1})$ with rational components,
which can be read off the   reduced Gr\"{o}bner basis of 
$\mfid^e =(Z_1-1, Z_2 -\frac{z_2}{z_1})$. 
One can check that $r(z_1,z_2)=r(1,\frac{z_2}{z_1})$ for
any $r\in \ratif=\K\left(\frac{z_2}{z_1}\right)$.  
\end{exy}

\begin{exy}\label{translation:algebraic} \textsc{translation+reflection}.
Consider the  group action from 
\ex{translation:def}, \ref{translation:graph},
\ref{translation:rational}, \ref{translation:cross}.
We chose the cross-section defined by  $\sid=(Z_1-Z_2)$ and found that
$\K(z_2^2)$ was the field of rational invariants.
Generic orbits have two components and the 
cross-section is of degree 2.  
Since
$\mfid^e=(Z_1-Z_2, Z_2^2-z_2^2)$, the two replacement invariants are 
$\xi^{(1)}=( z_2,z_2)$ and $\xi^{(2)}=(-z_2,-z_2)$.
Though rational functions, their components are not rational
invariants but only algebraic invariants.
Also $I^e=(Z_1-z_2,Z_2-z_2)\cap(Z_1+z_2,Z_2+z_2)$ is a reducible ideal of $\K(z)[Z]$, while
$I^G$ is an irreducible ideal of $\ratif[Z]$.
\end{exy}

\begin{exy}\label{rotation:algebraic} \textsc{rotation}.
Consider the  group action from \ex{rotation:def}, 
\ref{rotation:graph}, \ref{rotation:rational}, \ref{rotation:cross}.
We chose the cross-section defined by  $\sid=(Z_2)$.
Here the cross-section is again of degree 2 but the generic orbits
have a single component.  
Since $\mfid^e=(Z_2, Z_1^2-z_1^2-z_2^2)$ 
the two replacement invariants associated to $\sva$ 
are $\xi^{(\pm)}=(0,\pm\sqrt{z_1^2+z_2^2})$. 
\end{exy}

\subsection{Invariantization}
\label{inv}

In this section we introduce invariantization as a projection from the
ring of univariate polynomials over $\K[z]$ to the ring of univariate
polynomials over $\ratif$. It depends on the choice of a cross-section
and is computable by algebraic elimination.
As this projection extends to univariate polynomials over $\K(\sva)$
it can be understood as the computable counterpart to the  isomorphism 
\( \overline{\K(\sva)}\cong\algif \) that follows from \rprop{d-ext}.

We assume throughout this section that the field $\K$ is of
characteristic zero.
The ideal of the cross-section $\sva$ is taken alternatively in
$\K[z]$ and in $\K[Z]$. To avoid confusion we shall use in this
section $\sid[z]$ and $\sid[Z]$ to distinguish the two cases.
The localization of $\K[z]$ at $\sid[z]$ is denoted by $\locs$. 
In the proof of \theo{rewrite2} we have shown  that
$\ratif\subset\locs$.

The first approach for invariantization that draws directly on \cite{olver99} 
is to consider a replacement invariant  $\ninv$  associated to $\sva$ and 
the following chain of  homomorphisms: 
\begin{equation}\label{icv} 
\begin{array}{ccccc}
    \locs &\stackrel{\pr} \longrightarrow & \K(\sva) &\stackrel{{\isv}} \longrightarrow & \algif \\
     r(z) & \longmapsto & r(z)+\sid[z] & \longmapsto & r(\ninv)
\end{array} 
\end{equation}
The restriction of $\inv_\ninv={\isv}\circ\pi\colon \locs\to\algif$  to 
$\ratif$ is the identity map by
\theo{replacement}.
 We call the image  of a  rational function  $r(z)\in\locs$ 
under  $\inv_\ninv$  its \emph{\ninv-invariantization}.

If $\sva$ is a cross-section of degree $d$ there are  
$d$ distinct  associated replacement invariants 
$\ninv^{(1)}, \ldots, \ninv^{(d)}$. 
The image $\inv_\ninv(r(z))=r(\ninv)$ depends on the chosen
replacement invariant $\ninv$. 
Such is not the case of  
the minimal polynomial of $r(\ninv)$ over $\ratif$ which  depends
only on $\sva$ as we shall see. 
We  therefore define the $\sva$-invariantization as a map
taking a univariate polynomial over $\locs$ to a univariate polynomial
over $\ratif$.  The connection to the smooth invariantization of
\cite{olver99} is developed in \sec{felsolver}.

\begin{defy}\label{inv:def} The \sva-invariantization $\inv\alpha$ of
  a monic univariate polynomial $\alpha\in \locs[\dv]$ 
 is the squarefree part of
$\prod_{i=1}^d\alpha(\ninv^{(i)},\dv)$, 
where $\ninv^{(1)}, \ldots, \ninv^{(d)}$ are  
the $d$ replacement invariants
associated to the cross-section $\sva$. 
\end{defy}

Readers familiar with computer algebra techniques  
can see that $\inv \alpha$ belongs to  $\ratif[\dv]$ with the
following line of arguments.
The replacement invariants $\ninv^{(1)}, \ldots, \ninv^{(d)}$ are the
$d$ distinct zeros  of the zero dimensional prime ideal
$\mfid^G$ of $\ratif[Z]$. 
By a transcription of the primitive element
theorem, see for instance \cite[Proposition 4.2.2-3]{greuel}, 
they are thus the images by a polynomial map  
$\psi: \theta \mapsto (\psi_1(\theta), \ldots,\psi_n(\theta))$  over $\ratif$ 
of the roots ${\theta}^{(1)}, \ldots, {\theta}^{(d)} \in \algif$ 
of an irreducible  univariate polynomial of   degree $d$ with 
coefficients in $\ratif$. The coefficients of the polynomial  
$$ \prod_{i=1}^d\alpha(\ninv^{(i)},\dv) =  \prod_{i=1}^d\alpha(\psi({\theta}^{(i)}),\dv)$$ 
are elements of the field extension
 $\ratif(\theta^{(1)}, \ldots,\theta^{(d)})$
of $\ratif$ that  are invariant under all permutations of the
$\theta^{(i)}$.
By \cite[Section 8.1]{waerden71} or \cite[Theorem 8.15]{geddes92},
that polynomial belongs to  $\ratif[\dv]$ and thus so does  
 its squarefree part $\inv\alpha$   \cite[Section 8.1]{waerden71}.

 For a Galois theory oriented reader the details are given below. 
By definition $\inv\alpha$ belongs to the extension $\K(\ninv^{(1)},
\ldots, \ninv^{(d)})$, which we denote by $\K_\xi$.  Due to
\theo{replacement}  $\K_\xi=\ratif(\ninv^{(1)}, \ldots,
\ninv^{(d)})$. In order to  prove that $\inv\alpha \in \ratif[\dv]$ we
will show  that this polynomial is preserved by the Galois group of
the extension $\K_\xi\supset \ratif$. We need the following
proposition.

\begin{propy} 
 Let $\{\ninv^{(1)},\ldots, \ninv^{(d)}\}$ be the set of replacement
 invariants corresponding to a cross-section $\sva$ of degree
 $d$. Then the field $\K_\ninv=\K( \ninv^{(1)},\ldots, \ninv^{(d)})$
 is a splitting field of a univariate polynomial $\beta(z,\dv)\in
 \ratif[\dv]$ of degree $d$.  The Galois group of the extension
 $\K_\ninv\supset \ratif$  permutes the $n$-tuples
 $\ninv^{(1)},\ldots, \ninv^{(d)}$.  
 \end{propy}

 \begin{proof}  
Due to the replacement Theorem~\ref{replacement} one has the equality
$\K(\ninv^{(1)})=\ratif(\ninv^{(1)})$. From Corollary~\ref{d-ext} it
follows that   $\ratif(\ninv^{(1)})$ is an extension of degree $d$ of
$\ratif$ for $i=1..d$. Since $\K$ assumed to be of characteristic
zero, the components  $\ninv^{(1)}_1,\dots,\ninv^{(1)}_n$  of n-tuple
$\ninv^{(1)}$ are separable over $\ratif$. Hence there exists a
primitive element $\prel_1\in\K(\ninv^{(1)})$, such that
$\K(\ninv^{(1)})=\ratif(\ninv^{(1)})=\ratif(\prel_1)$, where $\prel_1$
is a root of  an irreducible univariate polynomial
$\beta(z,\dv)\in\ratif[\dv]$ of degree $d$  \cite[Theorem 5.4.1]{cox04}. 

Let $\sigma_{ji}\colon \K(\ninv^{(i)})\to\K(\ninv^{(j)})$ be the
$\ratif$-isomorphism induced by exchanging $\ninv^{(i)}$ and
$\ninv^{(j)}$. Then  $\prel_j=\sigma_{j1}(\prel_1)$ is a primitive
element of the extension $\K(\ninv^{(j)})\supset\ratif$. Indeed, since
$\prel_1$ is the primitive element of $\ratif(\ninv^{(1)})$,   for each
$i=1..n$, there exists  polynomial   $\psi_i$ over $\ratif$ such that
$\ninv^{(1)}_i=\psi_i(\prel_1)$.  
Since $\sigma_{j1}$ is a $\ratif$-isomorphism, it follows that 
$\ninv^{(j)}_i=\sigma_{j1}(\ninv^{(1)}_i)=\sigma_{j1}(\psi_i(\prel_1))=\psi_i(\sigma_{j1}(\prel_1))=\psi_i(\prel_j)$
for $i=1..n$.   Thus $\prel_j$ is a primitive element of
$\K(\ninv^{(j)})\supset\ratif$, and so
 $\K_\ninv=\ratif( \prel_1,\dots,\prel_d)$  

In addition, we proved that $n$-tuples  $\ninv^{(1)},\ldots,
\ninv^{(d)}$ are images of  $\prel_1,\dots,\prel_d$  under the
polynomial map
$\psi=(\psi_1,\dots\psi_n):\algif\to\left[\algif\right]^n$, where the
coefficients  of  the univariate polynomials $\psi_1,\dots\psi_n$ are
in $\ratif$. 
Since    $\ninv^{(1)},\ldots, \ninv^{(d)}$ are distinct tuples, then
$\prel_{1},\ldots, \prel_{d}$ are distinct elements of $\algif$. 
We will now show that $\prel_1,\dots,\prel_d$ are roots  of the minimal
polynomial $\beta\in\ratif[\dv]$ that defines $\prel_1$. 

Indeed, since the field $\ratif$ is fixed under $\sigma_{j1}$, for
$j=1..d$, then so is the polynomial $\beta$. Thus
$\prel_j=\sigma_{j1}(\prel_1)$ are  roots of the polynomial $\beta$. It
follows that $\K_\ninv=\ratif(\prel_1,\dots,\prel_d)$  is the splitting
field of an irreducible univariate polynomial $\beta\in\ratif[\dv]$ of
degree $d$.

The elements of the $Gal(\K_\xi/\ratif)$  permute the roots
$\prel_1,\dots,\prel_d$ of the polynomial $\beta$, and therefore it
permutes the tuples  $\ninv^{(j)}=\psi(\theta_j)$ for all $j=1..d$.  
\end{proof} 

\begin{coly} Let $\alpha(z,\dv) \in\locs$ be a univariate polynomial
  over $\locs$. Then its $\sva$-invariantization $\inv\alpha$ is a
  polynomial over $\ratif$.  
\end{coly}

\begin{proof} The Galois group of the extension $\K_\ninv\supset
  \ratif$ induces permutations of the  n-tuples $\ninv^{(1)},\ldots,
  \ninv^{(d)}$. 
Thus the polynomial
$p(\dv)=\prod_{i=1}^d\alpha(\ninv^{(i)},\dv)\in\K_\ninv[\dv]$ is fixed
under  $Gal (\K_\ninv/\ratif)$. Hence its coefficients belong to
$\ratif$. By definition $\inv\alpha$ is the square-free part of
$p(\dv)$, and hence it  is also fixed under the Galois group,  
since it has the same roots in $\K_\ninv $ as $p(\dv)$ itself
\cite[Proposition 5.3.8]{cox04}, and the Galois group permutes these
roots. Thus  its coefficients of $\inv\alpha$ are in $\ratif$.  
\end{proof}

The  following  properties follow directly from the definition of the map $\inv$:
\begin{enumerate}
\item A $\algif$-zero of $\inv\beta$ 
      is a $\algif$-zero of a $\beta(\xi^{(i)}, \dv)$ and conversely. 
\item If $\beta \in \ratif[\dv]$ then $\inv\beta=\beta$ since 
   $\beta(\ninv^{(i)},\dv)=\beta(z,\dv)$  by \theo{replacement}.
\item If $\alpha\equiv\beta\mod\sid[z]$ then $\inv\alpha=\inv\beta$ since
      the elements of $\sid[z]$ vanish on all $\ninv^{(i)}$. 
\end{enumerate}

The last property shows that $\inv$  induces a map $\phi$ 
from the set of monic polynomials of $\K(\sva)[\dv]$ 
to the set monic polynomials of $\ratif[\dv]$ s.t. $\inv =  \phi\circ\pi $.

From the first property it follows that 
$\beta(\ninv^{(i)},\dv)$ divides $\inv\beta(z,\dv)$ in
$\K(\ninv^{(i)})[\dv] \supset \ratif[\dv]$ when 
$\beta(\ninv^{(i)},\dv)$ is squarefree. 
Since $\K(\sva)\cong\K(\ninv^{(i)})$ this amounts to the following proposition that 
is used in \sec{felsolver}.

\begin{propy}\label{factor} 
 Let $\beta$ be a monic polynomial of $\locs[\dv]$. 
 If  $\beta$ is squarefree
 when considered in $\K(\sva)[\dv]$
 then it divides $\inv\beta(z,\dv)$ in   $\K(\sva)[\dv]$, 
that is 
 there exists $q(z,\dv) \in\locs[\dv]$ 
 s.t. $\inv\beta(z,\dv)\equiv q(z,\dv)\beta(z,\dv) \mod  \sid[z]$.
\end{propy}

Also we recognize in the definition of the invariantization map the
norm of a polynomial in a algebraic extension   
\cite[Section 8.8]{geddes92}.
We reformulate the results extending those of that text namely:
\begin{itemize}
\item[-] $\inv\beta$ can be computed by algebraic elimination. 
\item[-] if $\beta(\ninv^{(i)},\dv)$ is the minimal polynomial over 
$\K(\ninv^{(i)})\subset\algif$ of an  element  in $\algif$, 
then $\inv\beta$ is the minimal polynomial of this element over
$\ratif$ 
\end{itemize}

The algebraic elimination to compute  $\inv\beta$  can be performed by  several techniques. 
For a strict generalization of  \cite[Section 8.8]{geddes92} one could
introduce a resultant formula, as developed in \cite{andrea05}.
We propose here  a formulation in terms of elimination ideals.

\begin{propy} \label{inv:comp} 
Let $\beta \in \locs [\dv]$ be a monic polynomial. Then 
its $\sva$-invariantization  $\inv\beta$
is the squarefree part of the 
monic generator of $(\mfid^G+\alpha(Z,\dv))\cap \ratif[\dv]$
where $\alpha(z,\dv)\in \K[z][\dv]$ is the numerator of $\beta$.
\end{propy}

\begin{proof} 
The leading coefficient of  $\alpha(Z,\dv)\in\K[Z][\dv]$  does not belong to $\sid[Z]$, and
therefore it does not belong to $\mfid^G$. It follows that
$(\mfid^G+\alpha(Z,\dv))\cap \ratif[\dv] \neq (0)$ since $\mfid^G$ is zero-dimensional.

Let $\gamma(z,\dv)$ be the monic generator of $(\mfid^G+\alpha(Z,\dv))\cap\ratif[\dv]$.
We first prove that $\inv\beta$ divides the squarefree part of
$\gamma(z,\dv)$. The fact that $\gamma(z,\dv)$ belongs to
$\mfid^G+\alpha(Z,\dv)$ can be 
written as $\gamma(z,\dv)\equiv q(z,Z,\dv)\alpha(Z,\dv) \mod \mfid^G$ where
$q(z,Z,\dv)\in \ratif[Z,\dv]$. Substituting $\ninv^{(i)}$ for $Z$ we have
$\gamma(z,\dv)=q'(z,\ninv^{(i)},\dv)\beta(\ninv^{(i)},\dv)$
where $q(z,\ninv^{(i)},\dv)$ and $q'(z,\ninv^{(i)},\dv)$ 
differ by the factor in $\K[\ninv^{(i)}]$ that 
distinguishes $\alpha(\ninv^{(i)},\dv)$ from $\beta(\ninv^{(i)},\dv)$.
Therefore all the factors $\beta(\ninv^{(i)},\dv)$ of $\inv\beta$  divide
$\gamma(z,\dv)$. Since  $\inv\beta$ is the  squarefree product  of
$\beta(\ninv^{(i)},\dv)$  it divides the squarefree part of
$\gamma(z,\dv)$. 

Conversely, we prove that the squarefree part of $\gamma(z,\dv)$ divides
$\inv\beta$. The $\algif$-zeros of $ \alpha (Z,\dv)+\mfid^G$ are the
$(n+1)$-tuples 
$(\ninv^{(i)},\brf_{i,j})$, where  $\brf_{i,j}$, $1\leq j\leq\deg\beta$, are the
roots of $\beta(\ninv^{(i)},\dv)$. Since  $\gamma(z,\dv)$ belongs to $
\alpha (Z,\dv)+\mfid^G$ its set of $\algif$-roots  
includes all the  $\brf_{i,j}$. Thus $\gamma$ and $\inv\beta$ have the same set of roots.
 Therefore the squarefree part of $\gamma$  divides $\inv\beta$
\end{proof}

Note that  the monic generator of $(\mfid^G+\alpha(Z,\dv))\cap
\ratif[\dv]$ is the monic generator of
$(\mfid^e+\alpha(Z,\dv))\cap\K(z)[\dv]$. 
 This latter is   an element of the reduced Gr\"{o}bner basis of
$\left( \alpha (Z,\dv)+\mfid^e \right)$ w.r.t a term order that
eliminates $Z$. It follows from  \rprop{cocorico} that
it belongs to $\ratif[\dv]$.
Therefore  computations over $\K(z)$  lead to the correct
reasult over  $\ratif$. 

The last proposition provides the computable counterpart of the
isomorphism \( \overline{\K(\sva)}\cong\algif \),
 elements of $ \overline{\K(\sva)}$ or $\algif$ being
represented by irreducible monic polynomials over $\K(\sva)$ or $\ratif$ respectively.

\begin{propy}\label {irr}  Let $\alpha$ be a monic polynomial of $\locs[\dv]$. 
 The polynomial $\inv\alpha\in \ratif[\dv]$ is irreducible if and only
 if $\alpha$ is a power of an irreducible polynomial  when considered in $\K(\sva)[\dv]$. 
\end{propy}

\begin{proof}
Note that $\inv(\beta\,\gamma)$,  
for $\beta,\,\gamma \in \locs[\dv]$, is the squarefree part of the product
$\inv \beta\,\inv \gamma$. So if $\alpha$ considered in 
$\K(\sva)[\dv]$ is the product of two relatively prime factors then 
$\inv \alpha$ cannot be irreducible.

We can replace $\alpha$ by its squarefree part when considered in
$\K(\sva)[\dv]$ without loss of generality and thus assume for the
converse implication that  $\alpha(z,\dv)$ is irreducible there.
Let $\bar{\alpha} \in\K[z][\dv]$ be obtained from $\alpha$ by cleaning
up the denominators. 
Then  $\bar{\alpha}(Z, \dv)$ is irreducible modulo $\mfid^G$
so that $\left(\bar{\alpha}(Z,\dv)+\mfid^G\right)$ is prime.
The monic generator  $\inv\alpha$ of
$\left(\alpha(Z,\dv)+\mfid^G\right)\cap\K(z)[\dv]$ is thus irreducible.
\end{proof}

The following example illustrates various properties of the $\sva$-invariantization map $\inv$.  

\begin{exy}  \textsc{scaling}. 
We consider the scaling action  defined in \ex{scaling:def} and 
  the cross-section defined by the ideal $\sid[Z]=(Z_1^2+Z_2^2-1)$. It is
  a cross-section of degree $2$. We have 
  $I^e=( Z_1^2-\frac{z_1^2}{z_1^2+z_2^2}, Z_2-\frac{z_2}{z_1}Z_1)$ and
  therefore the two replacement invariants are 
$$\ninv^{(\pm)} = 
  \left(\frac{\pm z_1}{\sqrt{z_1^2+z_2^2}},
  \frac{\pm z_2}{\sqrt{z_1^2+z_2^2}}\right). $$

The invariantization of $\alpha = \dv-z_{1}$ is 
$\inv\alpha= \dv^2-\frac{z_{1}^2}{z_{1}^2+z_{2}^2}$.
We have $ \inv\alpha = (\dv+z_1) \alpha +
\frac{z_1^2}{z_1^2+z_2^2}(z_1^2+z_2^2-1)\equiv  (\dv+z_1) \alpha \mod \sid[z]$.
We obtained $\inv\alpha$ by computing the reduced  Gr\"obner basis of the ideal
$(\dv-Z_{1}, Z_1^2-\frac{z_1^2}{z_1^2+z_2^2}, 
  Z_2-\frac{z_2}{z_1}Z_1 )$
with a term
order that  eliminates $Z_{1}$ and $Z_2$. Note that, although $\alpha$
defines a polynomial function, its invariantization defines two
algebraic invariants $\pm\frac{z_1}{\sqrt{z_1^2+z_2^2}}$. 

The invariantization of 
  $\beta= \dv^3+ \dv^2+z_2\dv+1$ is 
  $\inv\beta= \dv^6+2\,\dv^5+\dv^4+2\,\dv^3+\frac{z_2^2+2z_1^2}{z_1^2+z_2^2}\,\dv^2+1$.
We have  $\inv\beta\equiv(\dv^3+\dv^2-z_2 \dv+1)\beta \mod \sid[z]$.

In the next two instances the 
 monic polynomial is equal modulo $\sid[z]$ to a polynomial in 
$\algif[\dv]$. 
As a consequence, the invariantization equals to the original polynomial modulo $\sid[z]$ 

The polynomial $\gamma=\dv-z_1^{2}$ is equal to its \sva-invariantization   
$\inv\gamma = \dv-\frac{z_{1}^2}{z_{1}^2+z_{2}^2} \equiv  \gamma \mod \sid[z]$. 

The irreducible polynomial 
$\delta= \dv^2-\frac{z_1^2+z_2^2-1}{z_2^2}\dv-\frac{z_1^2}{z_2^2}$ 
becomes a reducible modulo $\sid[z]$:  
$\delta \equiv  \dv^2-\frac{z_1^2}{z_2^2}\mod \sid[z]$.
Its invariantization is thus reducible: 
\(\inv\delta =(\dv-\frac{z_1}{z_2})(\dv+\frac{z_1}{z_2}) 
 \equiv  \delta \mod \sid[z]\). 

\end{exy}

\section{Local invariants and the moving frame construction}
\label{felsolver}

In this section we connect the algebraic algorithms  presented in this
paper with their original  source of inspiration, the Fels-Olver
moving frame construction 
\cite{olver99}. It is shown in \cite{olver99} that in the case of a  \emph{ locally
  free smooth action} of a Lie group $\gva$ on a manifold $\zva$, a
choice of local cross-section  corresponds to a local
$\gva$-equivariant map $\rho\colon \zva\to \gva$. 
   This map provides a generalization of the classical geometrical
   moving frames\footnote{For this reason the map $\rho$ is called
     \emph {moving frame } in \cite{olver99}. We adopt the term  \emph
     {a moving frame  
map.}} \cite{cartan}. A moving frame map gives rise to an \emph{ invariantization}
process, a projection from the set of smooth functions to the set of local invariants.

We introduce an alternative definition of the smooth invariantization
process which, on one hand, generalizes the definition given in
\cite{olver99} to non-free, semi-regular actions and, on the other
hand, can be effectively reformulated in the algebraic context. We
make explicit comparisons with both the moving frame and the algebraic
constructions in  Section~\ref{smf} and Section~\ref{smooth::alg}
respectively.

In this section we consider real smooth manifolds. 
 All statements  and constructions from this section are applicable to
complex  manifolds.  In the latter case all maps and functions are assumed to be meromorphic.

\subsection{Local action of a Lie group on a smooth manifold} \label{LGaction}
We consider a Lie group $\gva$, with
identity denoted $e$, and a smooth manifold $\zva$ of dimension
$\dimz$.
We first review the necessary facts and terminology from the theory of Lie
group actions  on smooth manifolds. Our presentations is based on
\cite{GOV93,olver:yellow}. 

\begin{defy} A local action of a Lie group $\gva$ on a smooth manifold
  $\zva$ is a smooth map $g\colon\Omega\to\zva$, where
  $\Omega\supset\{e\}\times\zva$ is an open subset of
  $\gva\times\zva$, and  the map $g$ satisfies the following two
  properties: 
\begin{enumerate}
\item $g(e,\bz)=\bz$, $ \forall \bz\in \zva$.
 \item \(g(\bm, g(\bl, z))=g(\bm\cdot\bl, z)\), for all $\bz\in\zva$ and
$\bl,\bm\in\gva$ s.~t.~$(\bl,\bz)$ and $(\bm\cdot\bl,\bz)$ are in
$\Omega$. 
\end{enumerate}
\end{defy}

The \emph{orbit} of $\bz\in\zva$ is the image $\ova_{\bz}$ of the smooth map
$g_{\bz}\colon\gva\mapsto\zva$ defined by $g_{\bz}(\bl)=g(\bl,\bz)$.
The domain of $g_{\bz}$ is an open subset of $\gva$ containing $e$. 

For every point $\bz\in\zva$ the differential $dg_{\bz}:T\gva|_e
\rightarrow T\zva|_{\bz}$  maps   
the tangent space of $\gva$ at $e$ to the  tangent space of $\zva$ at
the point $\bz$. The tangent space  $T\gva|_e$ can be identified with
the Lie algebra ${\frak g}$ of $\gva$. Let ${\hat v} \in {\frak g}$ 
then  ${v}(\bz)= dg_{\bz}({\hat v})$ is a smooth vector field on $\zva$, called the  {\em 
infinitesimal generator} of  the  $\gva$-action corresponding to
${\hat v}$. The set of all infinitesimal generators for a
$\gva$-action form a Lie algebra, such that the map $\hat v\to v$ is a
Lie algebra homomorphism. 

By $\exp(\epsilon v,\bz)\colon \halala{R}\times\zva\to\zva$ we denote
the {flow} of $v$. The flow  is defined as an integral curve of the
vector field $v$ with the initial condition $\bz$.  One can prove that
every point of  the connected component of the orbit $\ova^0_{\bz}\ni
\bz$ can be reached from $\bz$ by a composition of flows of a finite
number of infinitesimal generators.

Let $\hat v_1,\dots,\hat v_\kappa$,
where $\kappa\geq\dimo$ is the dimension of the group, be a basis of
the Lie algebra of $\gva$. Then the  infinitesimal  generators
$v_1,\dots,v_\kappa$ span the  tangent space to the orbits at each
point of $\zva$.

\begin{defy} An action of a  Lie group $\gva$ on a smooth manifold
  $\zva$ is semi-regular if all orbits  have the same dimension.
\end{defy}

Throughout this  section the action
is assumed  to be semi-regular.  The
dimension of the orbits is denoted by $\dimo$.

\subsection{Local invariants}

We give definitions of local  invariants and  fundamental sets of
those.   We prove that the  existence of a fundamental set of local
invariants follows from the existence of a flat coordinate system. The proof is based on  standard arguments from differential geometry. 

\begin{defy}\label{linv} A smooth function $f$, defined on an open
subset $\ncs\subset\zva$, is a { local invariant} if $v(f)=0$ for any
infinitesimal generator $v$ of the $\gva$-action on
$\ncs$.
\end{defy}

Equivalently $f(\exp (\varepsilon v, \bz))=f(\bz)$ for all $\bz\in
\ncs$, all  infinitesimal generator $v$, and all real $\varepsilon$
sufficiently close to zero. 
If the group  $\gva$ is  connected, the function $f$ is
 continuous  on $\zva$,  and the condition of Definition~\ref{linv} is
 satisfied at every point of $\zva$ then $f$ is a 
 global invariant  on $\zva$ due to
 \cite[Proposition~2.6]{olver:yellow}. 
 In what follows we neither assume $f$ to be continuous outside of $\ncs$, nor $\gva$ to be
 connected. 
 
A collection of smooth functions $f_1,\dots,f_l$ are functionally
\emph{dependent} on a manifold $\zva$ if for each point $\bz\in\ncs$ there exists on open neighborhood  an $\ncs$  and   a non-zero differentiable  function $F$ in $l$ variables such
that $F(f_1,\dots,f_l)=0$ on $\ncs$. From the implicit function
theorem it follows that  $f_1,\dots,f_l$ are functionally {
  dependent} on $\ncs$ if and only if the  rank of the corresponding Jacobian
matrix  is less than $l$ at each point of $\zva$. We say that
functions   $f_1,\dots,f_l$ are \emph{independent} on $\zva$ if  they are not dependent when restricted  to any open subset of $\zva$.  As it is commented in \cite[p85]{olver:yellow} functional dependence and functional independence on $\zva$ do not exhaust the range of possibilities, except for analytic functions.  Throughout the section the term {\it independent functions}
means {\it functionally independent functions}. Finally we say that $f_1,\dots,f_l$  are independent at a point $\bz\in\zva$ if the rank of the corresponding Jacobian matrix is maximal at $\bz$.   Independence at $\bz$ implies independence on some open neighborhood of this point.    If $\ncs$ is an open subset of $\zva$ and  $f_1,\dots,f_n$ are independent at each point of $\zva$, then  these functions  provide a coordinate system on $\ncs$.

\begin{defy} A collection  of local invariants on $\ncs$ forms
 a \emph{fundamental set} if they are functionally independent, and any 
local invariant on $\ncs$ can be expressed as a smooth function of the invariants from this set.
\end{defy}

The Lie algebra of infinitesimal generators provides an integrable
distribution\footnote{An integrable distribution is a collection of
  smooth vector fields, whose span over the ring of smooth functions
  is closed with respect to Lie bracket.} of smooth vector-fields on
$\zva$, whose integral manifolds are orbits. For a semi-regular action
this distribution is of  constant rank $\dimo$, the dimension of the
orbits.  It follows from  Frobenius  theorem that in an open
neighborhood $\ncs$ of each point  there exists a coordinate system
$x_1,\dots, x_\dimo,y_1,\dots,y_{n-\dimo} $  
 such that the connected components of the orbits on $\ncs$
 are level sets of the last $n-\dimo$ coordinates  \cite[p.~262]{spivak70}
and \cite[Theorem~1.43]{olver:yellow}.
Such coordinate system is called {\em flat, or straightening}.  The
proof of the following theorem establishes that 
$y_1,\dots,y_{n-\dimo}$ form a  fundamental set of local
invariants.

\begin{theoy}\label{fund}
Let $\gva$ be a Lie group acting semi-regularly on an $n$-dimensional
manifold $\zva$. Let $\dimo$ be the dimension of the orbits. 
In the neighborhood of each point $\bz\in \zva$ 
 there exists a fundamental set of  $n-\dimo$ local  invariants.
\end{theoy}

\begin {proof} By  Frobenius theorem there exists a flat coordinate
  system  $x_1,\dots, x_\dimo$, $y_1,\dots,y_{n-\dimo} $  
 in a neighborhood $\ncs\ni \bz$.  
 The
  connected components of the orbits on $\ncs$ coincide  with the
  level sets of the last $n-\dimo$ coordinate functions. Thus 
  $y_1,\dots,y_{n-\dimo}$ are constant on  the connected
  components of the orbits, and therefore they are local invariants,
  being smooth   and functionally independent by  definition of a coordinate
  system. It remains to show that any other invariant is locally
  expressible in terms of them. Let $v$ be an infinitesimal generator of
  the
group action. Since $v(y_i)=0$ for $i=1..(n-\dimo)$ then
$v=\sum_{i=1}^{\dimo}v(x_i)\bv x i $ is a linear combination of the first
$\dimo$ basis vector fields.
Let $v_1=\sum_{i=1}^{\dimo}a_{1i}\bv x i ,\dots,
v_\kappa=\sum_{i=1}^{\dimo}a_{\kappa i}\bv x i $ be a basis of
infinitesimal generators of the group action. Without loss of
generality we may assume that the first $\dimo$ generators
$v_1,\dots, v_\dimo$ are linearly independent at each point of $\ncs$.
Let $f( x_1\dots, x_\dimo,y_1,\dots y_{n-\dimo})$ be a local
invariant, then 
$v_j(f)=\sum_{i=1}^r a_{ji}\frac{\partial f}{\partial  x_i } =0$ 
for $j=1..\dimo$.   This is a homogeneous system of $\dimo$
linear equation with $\dimo$ unknowns $ \frac{\partial f}{\partial
  x_1}, \dots,  \frac{\partial f}{\partial x_\dimo}$.  Since $v_1,\dots,
v_\dimo$ are linearly independent at each point, the rank of the
system is maximal. Thus    $ \left(\frac{\partial f}{\partial x_1}=0,
  \dots,  \frac{\partial f}{\partial x_n}=0\right)$ is the only
solution. Hence $f$ is a function of invariants $y_1,\dots, y_{n-\dimo}$.
 \end{proof}

 The  existence of  a fundamental set of local invariants, therefore, follows from the existence of a flat
coordinate system. The proof is not constructive however.  The invariantization process,
introduced in next section, leads to a different characterization of a
fundamental set of invariants. Invariantization, and therefore fundamental invariants, can be effectively computed 
either by the algorithms   of \sec{invariantization}, in the case of a 
rational action of an algebraic group (see \sec{smooth::alg}),
or by  the moving frame method  of \cite{olver99}, in the case of 
a locally free action  of a Lie group (see \sec{smf}).

\subsection{Local cross-section and smooth invariantization}\label{smooth::section}

We  define local cross-sections to the orbits  and show that a local
cross-section  passing through any given point  can  easily be constructed.
A local cross-section gives rise to an equivalence relationship
on the ring of smooth functions such that any class has a single 
representative that is a local invariant.
This leads to an \emph{invariantization map},  a projection from the ring of smooth functions
to the ring of local invariants. 
It generalizes the invariantization process defined in
\cite{olver99} to semi-regular actions.  
 Although a possibility of such generalization is indicated
 in the remarks of  \cite[Section~4]{olver99}, the precise  definitions and theorems, appearing in this section, are new.

\begin{defy} \label{lsection}
An embedded submanifold $\sva$ of $\zva$ is a \emph{local cross-section}
 to the orbits if there is an open set $\ncs$ of $\zva$
 such that
\begin{itemize}
\item[-] $\sva$ intersects $\ova_{\bz}^0\cap\ncs$ at a unique
      point $\forall\bz\in\ncs$, where $\ova_{\bz}^0$ is the connected component of  $\ova_{\bz}\cap
      \ncs$, containing $\bz$.
\item[-] for all $\bz \in \sva\cap \ncs$,
   $\ova_{\bz}^0$ and $\sva$ are transversal and of complementary
   dimensions.
\end{itemize}
\end{defy}

The second condition in the above definition is
equivalent to the following condition on tangent spaces:
$T_{\bz}\zva=T_{\bz}\sva\oplus T_{\bz}\ova_{\bz}$, $\forall \bz\in
\sva\cap \ncs$.

An embedded submanifold  of codimension  $\dimo$ is locally 
given as the zero set of  $\dimo$ independent functions. Assume that 
 $h_1(z),\dots, h_\dimo(z)$ define $\sva$  on $\ncs$.
The tangent
space at a point of $\sva$ is  the kernel of the Jacobian
matrix $J_h$ at this point.  A basis of infinitesimal  generators
$v_1,\dots,v_\kappa$,  where $\kappa\geq\dimo$ is the dimension of the
group, span the  tangent space to the orbits at each point of
$\sva$. Therefore 
the submanifold $\sva$ is a local cross-section if and only if the
span of the infinitesimal generators $v_1,\dots,v_\kappa$ 
has a trivial  intersection with the kernel of $J_h$ on $\sva$. 
Equivalently: 
 \begin{equation}\label{tr} \mbox{ the rank of the } \dimo\times\kappa
   \mbox{ matrix  }
   \left(v_j(h_i)\right)_{i=1..\dimo}^{j=1..\kappa}=J_h \cdot V 
   \mbox{ equals to } \dimo \mbox{ on } \sva,\end{equation}
  where $V$ is the $n\times \kappa$ matrix, whose $i$-th column
  consists of the coefficients of the  infinitesimal generator $v_i$
  in a local coordinate system.
In the next theorem we prove the existence of a local cross-section through every point. The first paragraph of the proof   provides a  {simple practical algorithm to construct a coordinate local cross-section through a point}.  An
algebraic counterpart of this statements is given by Theorem~\ref{cross:existence}.
 
\begin{theoy} \label{banana} Let $\gva$ act semi-regularly on $\zva$.
 Through every point $\bz\in\zva$ there is
  a local cross-section  that is defined as the level set of $\dimo$
  coordinate functions.
\end{theoy}

\begin{proof}
Let $V$ be the $n\times \kappa$ matrix of the coefficients of
the infinitesimal generators $v_1,\dots,v_\kappa$
relative to a coordinate system $z_1,\dots,z_k$.  The rank of $V$ equals to the dimension of the orbits $\dimo$.
Thus there exist $\dimo$ rows of $V$ that form an $\dimo\times\kappa$ submatrix $\hat V$ of
rank $\dimo$ at the point $\bz$, and therefore it has rank $\dimo$ on
an open neighborhood $\ncs_1\ni \bz$. Assume that these rows
correspond to coordinate $z_{i_1},\dots,z_{i_\dimo}$. Let
$(c_1,\dots,c_n)$ be coordinates of the point $\bz$, then functions 
$h_1=z_{i_1}-c_{i_1},\dots,h_\dimo=z_{i_\dimo}-c_{i_\dimo}$ satisfy
condition \Ref{tr}.  The common zero set $\sva$ of these functions
contains $\bz$. 

It remains to  prove that there  exists a neighborhood $\ncs\ni\bz$
such that $\sva$ intersects each connected component  of the orbits on
$\ncs$ at a unique point.  Let  $x_1,\dots,x_s,y_1,\dots,y_{n-\dimo}$ be a flat coordinate system  in an open
neighborhood $\ncs_2\ni \bz$. Due to \theo{fund}  $y_1,\dots,y_{n-\dimo}$ are independent local invariants.  We
will show that  functions
$z_{i_1},\dots,z_{i_\dimo},y_1,\dots,y_{n-\dimo}$  provide  a  coordinate system  an open set
$\ncs=\ncs_1\cap\ncs_2$ containing $\bz$.  Without loss of generality
we may assume that $\{z_{i_1},\dots,z_{i_\dimo}\}=
\{z_{1},\dots,z_{\dimo}\}$ are the first $\dimo$ coordinates. In terms of flat coordinates $z_i=F_i(x,y),\,i=1..\dimo$, where $F_i$ are smooth functions on $\ncs_2$. 
Since $v_i(y_j)=0$ for $i=1..\kappa,\,j=1.. n-\dimo$, then  \begin{equation}\label{au1}\left(v_j(z_i)\right)^{j=1..\kappa}_{i=1..\dimo}=\left(\frac{\partial F_i}{\partial x_r}\right)^{r=1..\dimo}_{i=1..\dimo}\cdot\left(v_j(x_r)\right)^{j=1..\kappa}_{r=1..\dimo}.\end{equation} 
  We note that $\left(v_i(z_j)\right)_{j=1..\dimo}^{i=1..\kappa}=\hat V$ is
 $\dimo\times\kappa$ matrix of rank $\dimo$ at each point of $\ncs$. Matrix   $\left(v_j(x_r)\right)^{j=1..\kappa}_{r=1..\dimo}$ also has maximal rank $\dimo$ on $\ncs$. Therefore the  matrix $\left(\frac{\partial F_i}{\partial x_r}\right)^{r=1..\dimo}_{i=1..\dimo}$ is invertible on $\ncs$.  By looking at the rank of the corresponding Jacobian matrix in flat  coordinates, we  conclude that functions $z_1,\dots,z_\dimo,y_1,\dots,y_{n-\dimo}$ are independent at each point of  $\ncs$, and therefore define a coordinate system on $\ncs$.

By construction all points on $\sva$ have the same
$z$-coordinates. Thus two distinct points of $\sva$ must differ by at
least one of the $y$-coordinates. Since 
$y$ coordinates are constant on the connected components of the orbits
on $\ncs$, distinct points of $\sva$   belong to distinct connected
components  of the orbits. 
\end{proof}

 Given a cross-section on $\ncs$ one can define a projection from the set of
smooth functions on $\ncs$ to the set of local invariants.

\begin{defy} \label{iota} Let $\sva$  be a local cross-section to
  the orbits  on an open set $\ncs$. Let $f$ be a smooth function on
$\ncs$. The \emph{invariantization} $\icv f$ of $f$  is the function
on $\ncs$  that is defined, for  $\bz\in \ncs$, by  $\icv f(\bz)=f (\bz_0),$ where
$\bz_0=\ova_{\bz}^0\cap\sva$. 
\end{defy}
In other words, the invariantization of a function $f$ is obtained by
spreading the  values of  $f$ on $\sva$ along the orbits.  The next
theorem shows that $\icv f$ is the unique local invariant with the
same  values on $\sva$ as $f$. 
\begin{theoy} \label{invariantizationII}
Let a Lie group $\gva$ act semi-regularly on a manifold $\zva$,
and let $\sva$ be a local cross-section.
 Then $\icv f$ is the
unique local invariant defined on $\ncs$
whose restriction to $\sva$ is equal to the restriction of $f$ to
$\sva$. In other words $\icv f|_\sva=f|_\sva$.
\end{theoy}

\begin{proof} For any  $\bz\in\ncs$ and small enough $\varepsilon$ the
  point $\exp(\varepsilon v, \bz)$ belongs to the same connected
  component $\ova_{\bz}^0$. Let 
$\bz_0=\ova_{\bz}^0\cap\sva$. Then $\icv f\left(\exp(\varepsilon v,
  \bz)\right)=f(\bz_0)=\icv f(\bz)$, and thus $\icv f$ is a local
invariant.  By definition  $\icv f(\bz_0)=f(\bz_0)$ for all $\bz_0\in
\sva$. 

In order to show its smoothness we write  $\icv f$ in terms of  flat
coordinates $ x_1,\dots, x_\dimo$, $y_1,\dots,y_{n-\dimo}$. 
By
 probably shrinking $\ncs$, we may assume that $\sva$ is given by  the
 zero-set of  smooth independent functions $h_1(x_1,\dots,
 x_\dimo,y_1,\dots,y_{n-\dimo}),\dots,$ $ h_\dimo(x_1,\dots,
 x_\dimo,y_1,\dots,y_{n-\dimo})$. From the  transversality  condition
 \Ref{tr} and local invariance of $y$'s,   it follows that the first
 $\dimo$ columns of the Jacobian matrix $J_h$ form a  submatrix of
 rank $\dimo$.
Thus the cross-section $\sva$ can be described   as a graph

$x_1=p_1(y_1,\dots,y_{n-\dimo}), \dots,
x_\dimo=p_\dimo(y_1,\dots,y_{n-\dimo})$, where 
$p_1,\dots,p_\dimo$ are smooth functions. Then the function
\begin{small} $$\icv f(x_1,\dots, x_\dimo,y_1,\dots,y_{n-\dimo} )=f\left(
   p_1(y_1,\dots,y_{n-\dimo}), \dots,
   p_\dimo(y_1,\dots,y_{n-\dimo}\right),y_1,\dots,y_{n-\dimo})$$
   \end{small} is
 smooth, as a composition of smooth functions. 

To
prove the uniqueness, assume that an invariant function $q$ has the
same values on $\sva$ as $f$, then the invariant function $h=\icv f-q$
has zero value on  $\sva$. A point 
$\bz\in\ncs$ can be reached from $\bz_0={\sva}\cap \ova_{\bz}^0$ by a
composition of flows defined by infinitesimal generators. Without loss
of generality, we may assume that it can be reached by a single flow
$\bz=\exp (\epsilon v, \bz_0)$, where 
$\exp (\varepsilon v, \bz_0)\subset \ova_{\bz}^0$ for all $0\leq
\varepsilon \leq \epsilon$. 
From the invariance of  $h$ it follows that  $h\left(\exp (\epsilon v,
  \bz_0)\right)=h(\bz_0)=0$. Thus $q(z)=\icv f(z)$ on $\ncs$. 
\end{proof} 

Theorem~\ref{invariantizationII} allows us to view the
invariantization process as a projection from the set of smooth
functions on $\ncs$ to the equivalence classes of functions with the
same value on $\sva$. Each equivalence class contains a  unique
local invariant. The algebraic counterpart  of this
point of view is described in Section~\ref{inv}.

The invariantization of differential forms can be defined in a similar
implicit manner. It has been shown in \cite{olver99,ko03} that
the essential information about the differential ring of invariants
and the structure of differential forms can be computed from the
infinitesimal generators of the action and the equations that define
the cross-section, without explicit formulas for invariants.

\subsection{Normalized and fundamental invariants} \label{normfund}

The \emph{normalized    invariants} introduced in  \cite{olver99} 
are the invariantizations of the coordinate functions. 
They have the replacement property.
 In the algebraic context they correspond to
 \emph{replacement} invariants defined in
Section~\ref{invariantization}. This correspondence is made precise
 by Proposition~\ref{smooth::replacement}. 
We show that a  set of normalized  invariants   contains a
fundamental set of local 
 invariants.     
 
 All results of this  subsection  are stated under the following
 assumptions. A manifold  $\sva$ is a local cross-section  
   to the $\dimo$-dimensional orbits of a
  semi-regular $\gva$-action on an open $\ncs\subset\zva$, and $\icv$ is the corresponding invariantization map.  The set $\ncs$ is a single coordinate chart on $\zva$ with coordinate functions $z_1,\dots,z_n$.
   By possibly shrinking $\ncs$ we may assume that $\sva$
  is  the zero set of $\dimo$ independent smooth functions.

 Since our definition of invariantization differs from \cite{olver99}
we restate and prove the replacement theorem.

 \begin{theoy} \label{thomas}
    If  $f(z_1,\dots,z_n)$ is  a
   local invariant on $\ncs$ then  
   $f(\icv z_1,\dots,\icv z_n)=f( z_1,\dots,z_n)$. 
 \end{theoy}

 \begin{proof}
  Since \phantom{o}  $\icv z_1|_\sva=z_1|_\sva,\dots,\icv z_n|_\sva=z_n|_\sva$, then
 $f(\icv z_1,\dots,\icv z_n)|_\sva=f( z_1,\dots,z_n)|_\sva$.  
Thus
 functions $f(\icv z_1,\dots,\icv z_n)$  and $f( z_1,\dots,z_n)$ are
 both local invariants and have the same value on $\sva$. By
 Theorem~\ref{invariantizationII} they coincide. 
 \end{proof}

 \begin{lemy} \label{syzygy0}
Let $\sva$ be a local cross-section on $\ncs$, given as the zero set
of  $s$ independent functions
 $h_1,\ldots,h_s$. 
Then
$h_1(\icv z_1, \ldots, \icv z_n)=0, \ldots, h_s(\icv z_1, \ldots, \icv z_n)=0$ on $\ncs$.  
If for a differentiable  $n$-variable function $f$ we have 
$f(\icv z_1, \ldots, \icv z_n)\equiv 0$ on an open subset of $\ncs$, 
then  there exits  open ${\cal W}\subset \ncs$  such that ${\cal W}\cap\sva\neq\emptyset$  and at each point of  ${\cal W}\cap{\sva}$   functions $f$, $h_1, \ldots, h_s$ are not independent.  \end{lemy}
 
 \begin{proof}
  Since $h(\icv z)|_\sva=\icv h(z)|_\sva$  and both functions are
  invariants, one has $h(\icv z)=\icv h(z)$ by
  \theo{invariantizationII}. The latter is zero since $h|_{\sva}=0$. 
  Assume now that there exits a differentiable function $f$  and an open subset of $\cal V\subset \ncs$ such that 
  $f(\icv z_1, \ldots, \icv z_n)\equiv 0$ on $\cal V$.
  Since $f(\icv z)=\icv f( z)$ is invariant,    there exists  an open ${\cal W}\supset {\cal V}$ such that $f(\icv z_1, \ldots, \icv z_n)\equiv 0$ on ${\cal W}$  and ${\cal W}\cap\sva\neq\emptyset$. We conclude that $f(z_1,\dots,z_n)\equiv 0$ on $\sva\cap {\cal W}$.
  In this case  $f$ cannot be  independent of $h_1,\ldots,h_s$ at any point of $\sva\cap W$ since
  otherwise this would imply that  $\sva$ is of dimension less then $n-s$.
\end{proof}

\begin{theoy}
Let $\sva$ be a local cross-section on $\ncs$, given as the zero set
of  $s$ independent functions.
  The set   $\{\icv z_1,\dots,\icv z_n\}$ of the invariantizations  
  of the coordinate functions $z_1,\dots,z_n$ contains 
  a fundamental set of $n-\dimo$ local invariants on  $\ncs$. 
\end{theoy}

\begin{proof}
Due to  the implicit function theorem,   after a possible shrinking $\ncs$ and  renumbering of the
  coordinate functions, we may assume that $\sva$ is the zero set of  the functions 
  $h_1(z)=z_1-p_1(z_{\dimo+1},\dots, z_\dimz), \ldots,
  h_\dimo(z)=z_\dimo- p_\dimo(z_{\dimo+1},\dots,z_\dimz)$. 
Therefore 
 $\icv z_1 = p_1(\icv z_{\dimo+1},\dots, \icv z_\dimz),\ldots,
  \icv z_\dimo = p_k(\icv z_{\dimo+1},\dots, \icv z_\dimz)$ 
by \theo{invariantizationII}.
From \theo{thomas} we can conclude that any local invariant can be
written in terms of $\icv z_{\dimo+1},\dots, \icv z_\dimz$.
Since for every differentiable non-zero   $n-s$-variable function $f$, functions $f(z_{\dimo+1},\dots,  z_\dimz), h_1(z),\dots,h_\dimo(z)$ are independent at every point of $\ncs$, then  by \rlem{syzygy0}, $\icv z_{\dimo+1},\dots, \icv z_\dimz$ are
functionally independent on $\ncs$.
\end{proof}
\subsection{Relation between the algebraic and the smooth  constructions}
\label{smooth::alg}
We establish a connection  between the smooth  and the algebraic constructions. 
We show that the normalized invariants (\sec{normfund}) can be viewed as 
smooth representatives of 
the replacement invariants (\sec{repl}), and  that  algebraic
invariantization (\sec{inv}) provides a constructive approach to
smooth invariantization  (\sec{smooth::section}).

To be at the intersection of the hypotheses of the smooth 
and the algebraic settings  we consider a real algebraic group,
that is the set of real points of an algebraic group 
defined\footnote{This implicitly means that we know the 
 ideal $G$ (\sec{agroupaction:def}) 
from a set of  generators with coefficients in $\Ri$.}
over $\Ri$. It is a real Lie group  
\cite[the Proposition in Chapter 3, Section 2.1.2]{springer89}. 
Lie groups appearing in applications often satisfy this property. 
We also assume that the 
local action is given by a rational map \Ref{action}, in
\sec{agroupaction:def}, that satisfies \hyp{groupaction:hyp}. This
guarantees  semi-regularity of the action on an open set $\zva$ of
$\Ri^n$ as 
the orbits of non-maximal dimension are contained in an algebraic set
defined by minors of the matrix $V$ of  \Ref{tr}, in \sec{smooth::section}.

In \sec{groupaction} to \ref{invariantization}  we assumed for convenience 
of writing that the field of coefficients $\K$ was algebraically closed.
Yet the  algebraic constructions of those sections 
require no extension of the field of definition of the group or the
action. With the initial  data described above, \theo{invgb} produces a set of rational invariants in
$\Ratif$ that generate  $\Ratif$ by \theo{rewrite}. 

Rational invariants are obviously local invariants. We show that so 
are smooth representatives of algebraic invariants. The following
definition formalizes the notion of a smooth representative of an
algebraic function.  


\begin{defy}  A smooth map $F:\ncs\rightarrow\Ri^k$ is a
  smooth zero of 
  $\{p_1, \ldots, p_\kappa\} \subset \Ri(z)[\zeta_1,\ldots, \zeta_k]$ 
  if the coefficients of the $p_i$ are well defined
  on $\ncs$ and $p_i(\bz,F(\bz))=0$ for all $\bz \in \ncs$.
  In this case we also say that $F$ is a smooth zero of the ideal
  $(p_1, \ldots, p_\kappa)$.
\end{defy}

\begin{propy}\label{smz}
Assume $F:\ncs\rightarrow\Ri^k$ is a
  smooth zero of $\{p_1, \ldots, p_\kappa\} \subset \Ratif[\zeta_1,\ldots, \zeta_k]$. 
If $(p_1, \ldots, p_\kappa)$ is a zero dimensional
  ideal then the components of $F$ are local invariants. 
\end{propy}

\begin{proof}
Let $p\in\Ratif[\zeta]$, that is
$p(z,\zeta)=\sum_{\alpha\in\Ni^n}a_{\alpha}(z)\zeta^{\alpha}$, where
$a_{\alpha}(z)\in \Ratif$. Assume that  $p(\bz,F(\bz))=0$  for all $\bz\in \ncs$. 
For any $\bz\in\ncs$ and an infinitesimal generator $v$
there exits $\epsilon>0$, such that $\eez\in\ncs$ whenever
$|\varepsilon|<\epsilon$. Then
$p(\eez,\smz(\eez))=\sum_{\alpha\in\Ni^n}a_\alpha(\eez)\smz(\eez)^\alpha=0$.
Since the coefficients $a_{\alpha}$ are invariant  
$\sum_{\alpha\in\Ni^n}a_{\alpha}(\bz)\smz(\eez)^\alpha= 0$ for all $\bz\in\ncs$ and
small enough $\varepsilon$. Thus for a fixed point $\bz$ all the values
$\smz(\eez)$ for all sufficiently small $\varepsilon$ are the common roots of the set of 
 polynomials $\{p_1, \ldots, p_\kappa\}$. Since by the assumption the number of roots is finite, we conclude that
$\smz(\eez)=\smz\left(\exp(0 v, {\bz})\right)=\smz( {\bz})$ and thus the components of  $\smz(z)$ are local invariants.
\end{proof}

It follows from  \theo{banana} that,  through every point of $\zva$,
 there exists a local cross-sections
defined by linear equations over $\Ri$. 
Conversely, we can  consider a cross-section  $\sva$,   defined 
over $\Ri$, that has non singular real points, meaning that
the real part has the same dimension as the complex part. 
 For  any point $\bz \in\zva\cap\sva$ where the rank of the matrix 
\Ref{tr} does not drop, there is a neighborhood $\ncs$ on which $\sva$ 
defines a local cross-section, and  such points are dense in $\sva$.

The  $\Ratif$-zero of the zero dimensional ideal 
$\mfid^G= ( G+P+(z-g(\lambda,z)) )\cap \Ratif [Z]$ are precisely the
replacement invariants. 
According to the  previous proposition the smooth zeros of this ideal are local
invariants. 
We characterize the tuple of  normalized invariants as one of them.

\begin{theoy}\label{smooth::replacement}   Let $\sva$ be an algebraic
  cross-section which, when   restricted to an open set $\ncs$, 
defines a smooth cross-section.   The tuple of normalized invariants  $\icv z=( \icv z_1,\dots,\icv z_n)$ is the
  smooth zero of the ideal $\mfid^G$ whose components agree with the
  coordinate functions on  $\sva\cap \ncs$.  
\end{theoy}

\begin{proof} Let $\bz\in \ncs$   be an arbitrary point, and let
  $\bz_0$ be the point  of  intersection of $\sva$ with the connected
  component of $\ova_{\bz}\cap\ncs$, containing $\bz$.  
  Then there exists $\bl$ 
  in the connected component of the identity of $\gva$, such that
  $\bz_0=\bl\bz$ so that  $(\bz,\bz_0)$ is a zero of the ideal $I=O+P$.    
  By definition $\icv z(\bz)=\bz_0$ and therefore $(\bz,\icv z(\bz))$ is a
  zero of the ideal $I$ for all $\bz\in\ncs$. 
  Equivalently   $\icv z$ is a smooth zero of $\mfid^G$. By
  \theo{invariantizationII}  it is  the unique tuple of local
  invariants that agree with the coordinate functions on $\sva\cap\ncs$.  
\end{proof}

Therefore a replacement invariant not only generates 
algebraic invariants but their smooth representatives also generate local invariants.

\begin{exy}\label{scaling:smooth} \textsc{scaling}. 
The  action defined in \ex{scaling:def} corresponds to the 
following action of the multiplicative group $\Ri^{*}$:
$$g : \begin{array}[t]{ccl} 
  \Ri^{*} \times \Ri^2 & \rightarrow &\Ri^2 \\ 
  (\lambda, z_1,z_2) & \mapsto & (\lambda z_1, \lambda z_2).
\end{array}$$
The action is semi-regular on $\Ri^2\setminus\{(0,0)\}$.

In  \ex{scaling:cross} we chose the cross-section $\sva$
defined by $z_1=1$. 
The cross-section being of degree 1 there is a single associated
replacement invariant that corresponds to the tuple 
$(1,\frac{z_2}{z_1})$ of rational invariants. 

Let $\ncs=\{(z_1,z_2) \in \Ri^2 \; |\;z_1\neq 0\}$.
The components of the smooth map $F : \ncs \rightarrow \Ri^2$ s.t. 
$F(z_1,z_2)=(1,\frac{z_2}{z_1})$ are the normalized invariants
for the local cross-section $\sva\cap \ncs$. 
\end{exy} 

\begin{exy}\label{translation:smooth} \textsc{translation+reflection}.
The  action defined in \ex{translation:def} corresponds to the 
following action of the Lie group $\Ri\times \{1,1\}$  given by
 $$ g : \begin{array}[t]{ccl} 
  \Ri \times \{1,1\} \times \Ri^2 & \rightarrow &\Ri^2 \\ 
  (\lambda_1,\lambda_2, z_1,z_2) & \mapsto & 
  \left( z_1 +\lambda_1, \lambda_2\,z_2 \right) .
  \end{array}
$$
The action is semi-regular on $\Ri^2$.

In  \ex{translation:cross} we  chose the cross-section $\sva$
defined by $z_2=z_1$. 
There are two replacement invariants associated to $\sva$:
$\xi^{(\pm)}=( \pm z_2,\pm z_2)$. 
They both correspond to 
smooth maps $F^{(\pm)}:\Ri^2 \rightarrow \Ri^2$ the components of
which are local invariants. 

Only $(z_2,z_2)$ coincides with the
coordinate functions on $\sva$, that  defines a 
local cross-section on $\ncs=\Ri^2$. 
The normalized invariants are thus $(z_2,z_2)$. 
\end{exy}

\begin{exy} \label{rotation:smooth} \textsc{rotation}
The  action defined in \ex{rotation:def} corresponds to the
following action of the additive group $\Ri$  given by
 $$g : \begin{array}[t]{ccl} 
  \Ri \times \Ri^2 & \rightarrow &\Ri^2 \\ 
  (t, z_1,z_2) & \mapsto & 
 \left(\frac{1-t^2}{1+t^2} z_1 -\frac{2t}{1+t^2} z_2 ,  
       \frac{2t}{1+t^2} z_1+\frac{1-t^2}{1+t^2} z_2\right) . 
  \end{array}
$$ 
The action is semi-regular on $\Ri^2\setminus\{(0,0)\}$.  
 
In  \ex{rotation:cross} we  chose the cross-section $\sva$
defined by $z_2=0$. 
The replacement invariants associated to the cross-section $\sva$ are
the $\overline{\Ri(z)}^G$-zeros of the ideal 
 $\mfid^G =(Z_2,Z_1^2-(z_1^2+z_2^2))$. 

The smooth  maps  $F^{(\pm)} : \Ri^2\setminus\{(0,0)\} \rightarrow  \Ri^2$ 
s.t. $F^{(\pm)}(z_1,z_2)=(0,\pm \sqrt{z_1^2+z_2^2})$ 
are smooth zeros of $I^G$. Their components are thus local invariants.

The cross-section $\sva$ defines a local cross-section for instance on
$\ncs=\Ri^2\setminus \{(z_1,z_2)  \; |\;z_1=0, z_2\leq0\}$. 
As $F^{(+)}|_{\sva\cap\ncs} = z_1$, the tuple of normalized invariants
are $(0, \sqrt{z_1^2+z_2^2})$ on $\ncs$.
\end{exy} 

We conclude this section by linking the smooth invariantization and 
the algebraic invariantization introduced in \sec{inv}.  Recall that
the algebraic invariantization was a map that associated a univariate
polynomial over $\Ratif$ to univariate polynomials over $\K[z]_\sva$ (Definition~ \ref{inv:def}).
 
\begin{theoy} Let $\sva$ be an algebraic cross-section which, when
  restricted to an open set $\ncs$,  defines a local cross-section. Let
  $f:\ncs\rightarrow \Ri$ be a  smooth zero of a 
 univariate  polynomial $\beta\in\K(z)[\zeta]$. The smooth
  invariantization $\icv f$  of $f$ 
is a smooth zero of the algebraic $\sva$-invariantization
  $\inv\beta\in\Ratif[\zeta]$ of $\beta$.  
\end{theoy}

\begin{proof} 
The polynomial $\inv\beta(z,\zeta)=\sum_{i=1}^{k}b_{i}(z)\zeta^{i}$, where
$b_{i}\in\ratif$. Any  
point $\bz\in\ncs$ can obtained from the point $\bz_{0}\in\sva$ by a
composition of flows along infinitesimal generators of the group
action.  The argument will not change if we assume that
 $ \bz=\exp(\varepsilon v,\bz_{0})$ is obtained by the flow along a
 single vector field. Then from the invariance of $b_{i}(z)$ and  
local invariance of $\icv\smr(z)$ it follows that  $\forall\bz\in\ncs$:
\begin{eqnarray*}
\inv\beta(\bz,\icv\smr(\bz))
&=&\sum_{i=1}^{k}b_i\left(\exp(\varepsilon v,\bz_{0})\right)f\left(\exp(\varepsilon
v,\bz_{0})\right)^i\\
  &=&\sum_{i=1}^{k}b_{i}(\bz_0)\icv\smr(\bz_0)^i
  =\inv\beta\left(\bz_0,\icv\smr(z_{0})\right), \,
 \mbox{ where }  \bz_0\in\sva\cap\ncs.
\end{eqnarray*}
 From Proposition~\ref{factor} it follows that $\inv\beta$ is
 divisible by $\beta$ when restricted to $\sva$.  Thus
 $\inv\beta({\bz}_0,\smr ({\bz}_0))=0,\quad \forall
 {\bz}_0\in\sva\cap\ncs$, since  $\beta(\bz,\smr (\bz))\equiv 0 $ on
 $\ncs$.  
It follows that  $\icv \smr(z)$ is a  smooth zero
of  a polynomial $\icv\beta(z,\zeta)\in \K(z)^G[\zeta]$. 
\end{proof}

In particular if $r(z)$ is a rational function that is well defined on
$\ncs$, then its smooth invariantization  $\icv r(z)$  is a smooth 
zero of the $\sva$-invariantization  $\inv (\zeta - r(z))$ of  the polynomial
$\zeta-r(z)$. To discriminate the right one we only need to check that
its value coincide with the one of $r(z)$ on $\sva\cap\ncs$.

\newpage
\subsection{Moving frame map}\label{smf}

We show that the invariantization map described in Section~\ref{inv}
generalizes the  invariantization process described in \cite{olver99}.
The latter is restricted  to  locally-free actions, and is based on
the existence  of  a local  $\gva$-equivariant map
$\rho\colon \ncs\to \gva$. 
Although  local freeness of the  action guarantees the existence of $\rho$, due to the implicit function theorem, it might not be explicitly
computable. We review the Fels-Olver construction, and prove that  in the case of locally free actions it is equivalent to the one presented in Section~\ref{smooth::section}.

\begin{defy}\label{lfree} An action of a Lie group $\gva$ on a manifold
$\zva$ is \emph{locally free} if for every point $\bz\in \zva$ its
isotropy group $\gva_{\bz}=\{\bl\in \gva|\bl\cdot\bz=\bz\}$ is
discrete.
 \end{defy}

Local freeness implies semi-regularity of the action, 
the dimension of each orbit being equal to the dimension of the
group.  
  Theorem~4.4 from \cite{olver99}, 
 can be restated as follows
in the case of locally free actions.

\begin{theoy}\label{lmf}
A Lie group $\gva$ acts  locally freely on $\zva$ if and only if
every point of $\zva$ has an open neighborhood $\ncs$ such that there exists a
 map $\rho\colon \ncs\to \gva$
that makes the following diagram commute.
Here the map $\bm\mapsto
\bm\cdot\bl^{-1}$ is chosen for the action of $\gva$ on itself, and
 $\bl$ is taken in a suitable neighborhood
(depending on the point of $\ncs$)
of the identity in $\gva$.
\[ \xymatrix{
   \ncs \ar[d]_{\rho}  \ar[r]^{\bl} & \ncs\ar[d]^{\rho}\\
   \gva \ar[r]_{{\bl}} & \gva }\]
\end{theoy}
The map
$\rho$ is locally $\gva$-equivariant, $\rho(\bl\cdot \bz)= \rho \cdot
\bl^{-1}$ for $\bl$ sufficiently close to the identity, and is called \emph{ a moving frame map}.
If $\sva$ is a cross-section, then the equation
   \begin{equation}\label{mfdef}
         \rho(\bz)\cdot\bz \in \sva,
    \end{equation}
uniquely defines $\rho(\bz)$ in a sufficiently small neighborhood of
the identity. In particular, $\rho(\bz_0)=e$ for all $\bz_0\in \sva$.
Reciprocally, a moving frame map defines a local cross-section
to the orbits: $\sva=\{\rho(\bz)\cdot \bz \;|\; \bz \in \ncs\}\subset \ncs$. 

In local coordinates, Condition  \Ref{mfdef} gives rise to {\em
  implicit equations} for expressing the group parameters in terms of
the coordinate functions on the 
manifold.  When the group acts locally freely, the local existence of
smooth solutions is guaranteed by the transversality condition and
the implicit function theorem.  Since the implicit function theorem is
not constructive,  we might nonetheless not be able 
to obtain explicit formulas for the solution.

In \cite[Definition~4.6]{olver99} the invariantization of a function
$f$ on $\ncs$ is defined as the function whose value at a point 
$\bz\in\ncs$ is equal to $f(\rho(\bz)\cdot \bz)$. 
Next proposition shows that this  moving frame  based definition
  of invariantization is equivalent to  Definition~\ref{iota} given in
  terms of cross-section. The advantage of the latter definition is that 
  it is not restricted  to locally free actions. 

\begin{propy}\label{foinv}
 Let  $\rho$ be a moving frame map on $\ncs$. Then
\(\icv f(\bz)=f(\rho(\bz)\cdot\bz).\) 
\end{propy}

\begin{proof}
Local  invariance of $f(\rho(z)\cdot
z)$ follows from the local equivariance of $\rho$, i.~e.  
for $\bl$ sufficiently close to the identity:  $$f\left(\rho(\bl\cdot \bz)\cdot(\bl\cdot 
\bz)\right)=f\left(\rho( \bz)\cdot\bl^{-1}\cdot(\bl\cdot 
\bz)\right)=f(\rho(\bz)\cdot
\bz.$$
 Since $\rho(z_0)=e$ then $ f(\rho(\bz_0)\cdot\bz_0)=f(\bz_0)$ for all  $\bz_0\in \sva$. Thus $f(\rho(z)\cdot
z)$ is locally invariant and equals to $f$, when restricted to $\sva$.
The conclusion follows from Theorem~\ref{invariantizationII}.
\end{proof}

Thus the moving frame map offers an approach to invariantization that
is constructive up to the resolution of the implicit equations given
by \Ref{mfdef}. 
In the algebraic case  
 the moving frame map is defined by the ideal
$$\mfmid^e= 
\left(\,\gid+\sid+(Z-g(\lambda,z))\,\right) \;\cap\;\Ri(z)[\lambda].$$
Indeed, if  $(\bz,\bl)$ is a zero of $M=M^e\cap\Ri[z,\lambda]$, in an appropriate open set of
$\zva\times\gva$, then $\bl\cdot{\bz}\in \sva$.
The action is locally free if and only if  $\mfmid^e$ is zero dimensional. 

In this case, the smooth zero $F:\ncs\rightarrow\gva$ of $M^e$,
that is the identity of the group when restricted to $\sva$,
provides a moving frame map $\rho$ on $\ncs$.

If one can obtain the map $\rho$ explicitly, the invariantization map can be computed using \rprop{foinv}.
Even  in this favorable case,  the expression for $\rho$ often involves
algebraic functions which can prove difficult to manipulate symbolically. 
The purely algebraic approach proposed in \sec{invariantization} 
is more suitable for symbolic computation.

\section{Additional examples} 
\label{examples}

We first consider  a linear action of $S\!L_2$ on $\K^7$ taken from 
\cite{derksen99}. That latter paper presents 
an algorithm to compute a set of generators
of the algebra of polynomial invariants for the linear action of a
reductive group.
The ideal $\oid = (\gid + (Z-g(\lambda,z))) \cap \K[z,Z]$, where now $g$ is a
polynomial map that is linear in $z$, is also central in the
construction as a set of generators of $\K[z]^G$ is obtained by
applying the Reynolds operator, which is a projection from $\K[z]$ to
$\K[z]^G$, to generators of $O+(Z_1,\ldots,Z_{\dimz})$, the ideal of
the null cone.

The fraction field of $\K[z]^G$ is included in $\ratif$ but does not need to 
be equal. Conversely there is no known algorithm to
compute $\K[z]^G=\ratif \cap \K[z]$ from the knowledge
of a set of generators of $\ratif$.

\begin{exy}  \label{derksen} We consider the linear action of
$S\!L_2$ on $\K^7$ given by the following polynomials of
$\K[\lambda_1,\ldots,\lambda_4,z_1,\ldots, z_7]$:
\[ \begin{array}{c} g_1 = \lambda_1 z_1+\lambda_2 z_2,\quad
g_2=\lambda_3z_1+\lambda_4z_2\; \\ g_3 = \lambda_1 z_3+\lambda_2
z_4,\quad g_4=\lambda_3z_3+\lambda_4z_4\; \\ g_5 =
\lambda_1^2z_5+2\lambda_1\lambda_2z_6+\lambda_2^2z_7, \\ g_6 =
\lambda_3\lambda_1 z_5+\lambda_1\lambda_4
+\lambda_2\lambda_3z_6+\lambda_2\lambda_4z_7,\; \\ g_7 =
\lambda_3^2z_5+2\lambda_3\lambda_4z_6+\lambda_4^2
\end{array}
\] 
the group being defined by
$\gid=(\lambda_1\lambda_4-\lambda_2\lambda_3-1) \subset
\K[\lambda_1,\lambda_2,\lambda_3,\lambda_4]$.

The cross-section defined by $\sid=(Z_1+1,Z_2,Z_3)$ is of degree one.
The reduced Gr\"{o}bner basis (for any term
order) of the ideal $\mfid^e \subset \K(z)[Z]$ is indeed given by
\( \{Z_1 + 1, Z_2, Z_3, Z_4-r_2, Z_5-{r_3},
Z_6-{r_4}, Z_7 - r_1\} \)
where 
\[ \begin{array}{c}
{r_1}={z}_{7}\,{{z}_{1}}^{2}-2\,{z}_{2}\,{z}_{6}\,{z}_{1}+{{
z_2}}^{2}{z}_{5}, 
\quad
{r}_{2}={z}_{3}\,{z}_{2}-{z}_{1}\,{z}_{4},
\\
\ds 
{r_3}={\frac {{{z}_{3}}^{2}{z}_{7}-2\,{z}_{6}\,{z}_{4}\,{z}_{3}
+{z}_{5}\,{{z}_{4}}^{2}}{ \left({z}_{1}\,{z_4} -{z}_{3}\,{z}_{2} \right) ^{2}}},
\quad 
{r}_{4}={\frac {{z}_{1}\,
{z_6}\,{z}_{4}-{z}_{1}\,{z}_{3}\,{z}_{7}+
{z_3}\,{z}_{2}\,{z}_{6}-{z}_{2}\,{z}_{5}\,{z}_{4}}{{z}_{1}\,{z}_{4}-{z}_{3}\,{z}_{2} }}
\end{array}
\]

By \theo{rewrite2}, \(\ratif=\K(r_1,r_2,r_3,r_4)\). In this case
the rewriting of any rational invariant in terms of 
$r_1,r_2,r_3,r_4$ consists simply of the 
substitution of $(z_1,z_2,z_3,z_4,z_5,z_6,z_7)$ by $(-1,0,0,r_2,r_3,r_4,r_1)$. 
We illustrate this by rewriting the five  generating  polynomial
invariants computed in \cite{derksen99} in terms of $r_1,r_2,r_3,r_4$:
\[ \begin{array}{c} 
{z_2}^{2}z_5-2\,z_2\,z_6\,z_1+z_7\,{z_1}^{2}=r_1,\quad
z_3\,z_2-z_1\,z_4=r_2,
\\ 
{z_3}^{2}z_7-2\,z_6\,z_4\,z_3+z_5\,{z_4}^{2}=r_3{r_2}^{2}, \,\,
z_1\,z_3\,z_7-z_3\,z_2\,z_6+z_2\,z_5\,z_4-z_1\,z_6\,z_4=r_4\,r_2,
\\
{z_6}^{2}-z_7\,z_5={r_4}^{2}-r_1\,r_3,\quad 
\end{array}
\]

The reduced Gr\"{o}bner basis of $\oid^e$, relative to the total
degree order with ties broken by reverse lexicographical order, has 9
elements: 

\[ \begin{array}{c} 
{{ Z_6}}^{2}
-{ Z_7}\,{ Z_5}+{ r_1}\,{ r_3}-{ r_4}^{2},
\quad
{ Z_6}\,{ Z_4}
+{ r_3}\,{ r_2}\,{ Z_2}-{ r_4}\,{ Z_4}-{ Z_3}\,{ Z_7},
\\
{ Z_5}\,{ Z_4}
-{ Z_3}\,{ Z_6}+{ r_3}\,{ r_2}\,{ Z_1}-{ r_4}\,{ Z_3},
\quad
{ Z_3}\,{ Z_2}
-{ Z_1}\,{ Z_4}-{ r_2},
\\
{ Z_2}\,{ Z_6}
-{ Z_1}\,{ Z_7}+{ r_4}\,{ Z_2}-{\frac{r_1}{ r_2}}\, Z_4,
\quad
{ Z_2}\,{ Z_5}
+{ Z_1}\,{ r_4}-{ Z_6}\,{ Z_1}-{\frac {r_1}{ r_2}}\,{ Z_3} ,
\\
{{ Z_2}}^{2}
+{\frac {{ r_1}}{{ r_3}\,{{    r_2}}^{2}}}\,{{ Z_4}}^{2}
-{\frac {{ Z_7}}{{ r_3}}}
-2\,{\frac {{ r_4}}{{ r_3}\,{ r_2}}}\,{ Z_4}\,{ Z_2},
\quad
{{ Z_1}}^{2}
-{\frac {{ Z_5}}{{ r_3}}}-
2\,{\frac {{ r_4}}{{ r_3}\,{ r_2}}}\,{ Z_3}\,{ Z_1}
+{\frac {{ r_1}}{{ r_3}\,{{ r_2}}^{2}}}\,{{ Z_3}}^{2}
\\
{ Z_2}\,{ Z_1}
-{\frac {{ r_4}}{{ r_3}}}-{\frac {{ Z_6}}{{ r_3}}}
+{\frac {{ r_1}}{{ r_3}\,{{ r_2}}^{2}}}\,{ Z_4}\,{ Z_3}
-2\,{\frac {{ r_4}}{{ r_3}\,{ r_2}}}\,{ Z_4}\,{ Z_1},
\end{array}
\]
Though this Gr\"obner basis is obtained without much difficulty,
the example illustrates  the advantage obtained by considering
the construction with a cross-section: $\mfid^e$ has a much simpler
reduced Gr\"obner basis than $\oid^e$. 
\end{exy} 

We finally take a classical example in differential geometry: the
Euclidean action on the second order jets of curves. The variables
$x,y_0,y_1,y_2$ stand for the independent variable, the
dependent variable, the first and the second derivatives respectively.
We shall recognize the  curvature  as the non constant component of a replacement invariant.

\begin{exy} We consider the group defined by
$\gid=(\alpha^2+\beta^2-1, \epsilon^2-1)
\subset\K[\alpha,\beta,a,b,\epsilon]$.
The neutral element is $(1,0,0,0,1)$, the
group operation is\\
$ (\alpha',\beta',a', b', \epsilon')\cdot
(\alpha,\beta,a, b,\epsilon) =
(\alpha\alpha'-\beta\beta', \beta\alpha'+\alpha\beta',
a +\alpha a'-\beta  b',
 b+\alpha a'+\alpha b', \epsilon\,\epsilon')$ and the inverse map 
$(\alpha,\beta,a, b)^{-1}=
(\alpha,-\beta, -\alpha \, a- b\beta, \beta \,a-\alpha b, \epsilon)$.  
The rational action on
$\K^4$ we consider is given by the rational functions:
\[ \begin{array}{c} \ds g_1 = \alpha x-\beta y_0 +a,
\quad g_2 = \epsilon \beta x+\epsilon \alpha y_0 + b,\\ 
\ds g_3 = \frac{\beta+\alpha y_1}{\alpha-\beta y_0} , 
\quad g_4 =\frac{y_2}{(\alpha-\beta y_0)^3}.
 \end{array} \]

We have 
\[ \oid = \left(\left( 1+ y_1^2 \right)^3 Y_2^2-\left(1+Y_1^2\right)^3 y_2^2\right)\] 
and if we consider the  the cross section defined by 
$\sid=(X,Y_0,Y_1)$  
the reduced Gr\"{o}bner basis of $\mfid^e = \oid^e+ \sid$ is
\[ \left\{X,Y_0,Y_1, Y_2^2-\frac{y_2^2}{(1+y_1^2)^3} \right\}.\] 
According to \theo{rewrite} or  \theo{rewrite2}
\[ \ratif=\K\left(\frac{y_2^2}{(1+y_1^2)^{3}}\right). \]

The two replacement invariants $\ninv=(\ninv_1, \ninv_2, \ninv_3, \ninv_4)$
associated to the cross-sections are given by
\[ \ninv_1=0,\ninv_2=0,\ninv_3=0,
\ninv_4=\pm \sqrt{\frac{y_2^2}{(1+y_1^2)^3}}. \]

\end{exy}

\nocite{algeom4}
\bibliographystyle{plain}

\end{document}